\documentclass[11pt,a4paper]{article}
\usepackage{ART-VEM-Enhanced}%
\usepackage{amsthm}
\usepackage[capitalize]{cleveref}
\renewcommand{\red}[1]{\textcolor{red}{#1}}
\renewcommand{\fix}[1]{#1}
\renewcommand{\ourfix}[1]{#1}
\renewcommand{\ourfixBis}[1]{#1}
\renewcommand{\fixBis}[1]{#1}
\crefname{equation}{}{}
\Crefname{equation}{}{}
\crefname{figure}{Figure}{Figures}
\author{Stefano Berrone, Andrea Borio, Francesca Marcon \footnote{The authors wish
    to thank professors Paolo Tilli, Alessandro Russo, Carlo Lovadina and Silvia Bertoluzza for
    their insightful comments and the anonymous reviewers for helping to improve
    the quality of the manuscript. The three authors are members of the
      Gruppo Nazionale Calcolo Scientifico (GNCS) at Istituto Nazionale di Alta
      Matematica (INdAM). The authors kindly acknowledge partial financial
      support by INdAM-GNCS Projects 2020, by the Italian Ministry of Education,
      University and Research (MIUR) through the ``Dipartimenti di Eccellenza''
      Programme (2018–2022) – Department of Mathematical Sciences
      ``G. L. Lagrange'', Politecnico di Torino (CUP:E11G18000350001) and
      through the PRIN 2017 project (No. 201744KLJL\_004).}}%
\title{Lowest order stabilization free Virtual Element Method for the 2D Poisson equation}%
\date{}

\newtheorem{theorem}{Theorem}
 \newtheorem{lemma}{Lemma}
 \newtheorem{remark}{Remark}
 \newtheorem{definition}{Definition}
 \newtheorem{proposition}{Proposition}

\begin{document}


\maketitle
\begin{abstract}
  We introduce and analyse the first order Enlarged Enhancement Virtual Element
  Method (E\textsuperscript{2}VEM) for the Poisson problem. The method allows
  the definition of bilinear forms that do not require a stabilization
  term\fix{, thanks to the exploitation of higher order polynomial projections
    that are made computable by suitably enlarging the enhancement (from which
    comes the prefix of the name E\textsuperscript{2}) property of local virtual
    spaces. The polynomial degree of local projections is chosen based on the
    number of vertices of each polygon}. We provide a proof of well-posedness
  and optimal order a priori error estimates. Numerical tests on convex and
  non-convex polygonal meshes confirm the \ourfix{criterium for well-posedness
    and the} theoretical convergence rates.
\end{abstract}
\section{Introduction}
In recent years, the study of polygonal methods for solving partial differential
equations has received a huge attention. The main reason for this great interest
relies in the flexibility of polygonal meshes to discretize domains with high
geometrical complexity.  A large number of families of polygonal/polyhedral
methods has been developed, among them we can list Discontinuous Galerkin
Methods \cite{DiPietro2012,Riviere2008,HesthavenJanS2008NdGm}, Polygonal Finite
Elements (PFEM) \cite{PFEM_2004}, Mimetic Finite Difference Methods (MFD)
\cite{MFD_Beirao_Manzini,MFD_Brezzi,MimeticFiniteDifference_1981}, Hybrid High
Order Methods (HHO) \cite{DIPIETRO2015b,DIPIETRO2015a,DIPIETRO2014}, Gradient
Discretisation Methods \cite{Droniou2013,Droniou2018}, CutFEM \cite{BURMAN2015},
other methods that help in circumventing geometrical complexities are Extended
FEMs (XFEM) \cite{XFEM_1999}, Generalised FEMs (GFEM)
\cite{STROUBOULIS200043,STROUBOULIS20014081,STROUBOULIS200047} as well as
Ficticious Domain Methods \cite{GLOWINSKI1994283,Babuska2005}, Immersed Boundary
Methods \cite{Peskin2002}, PDE-constrained Optimization Methods
\cite{BPSa,BPSb,BPSc} and many others.  One of the most recent developments in
this field is the family of the Virtual Element Methods (VEM). These methods
were first introduced in primal conforming form in \cite{Beirao2013a} and were
later on applied to most of the relevant problems of interest in applications,
such as advection-diffusion-reaction equations
\cite{Beirao2015b,BBBPSsupg,BBM}, elastic and inelastic problems
\cite{Beirao2015a}, plate bending problems \cite{Brezzi2013a}, parabolic and
hyperbolic problems \cite{Vacca2015,Vacca2017}, \ourfix{simulations on unbounded domains \cite{DESIDERIO2021296}, }simulations in fractured media
\cite{BBS,BBBPS,BBBchapter}.

Standard VEM discrete bilinear forms are the sum of a singular part maintaining
consistency on polynomials and a stabilizing form enforcing coercivity. In the
literature, the stabilization term has been extensively studied, for instance in
\cite{Beirao2017}, and remains a somehow arbitrarily chosen component of the
method with several possible effects on the stability and conditioning of the
method.  Moreover, the stabilization term causes issues in many theoretical
contexts. The first one that we mention is the derivation of a posteriori error
estimates \cite{Cangiani2017,BBapost}, where the stabilization term is always at
the right-hand side when bounding the error in terms of the error estimator,
both from above and from below. Moreover, the isotropic nature of the
stabilization term becomes an issue when devising SUPG stabilizations
\cite{BBBPSsupg,BBM}, in problems with anisotropic coefficients, or in the
derivation of anisotropic a posteriori error estimates \cite{ABBDVW}. Finally,
other contexts in which the stabilization may induce problems are multigrid
analysis \cite{Multigrid_Antonietti} and complex non-linear problems
\cite{Wriggers2019}.

In this work, we introduce a new family of VEM, that we call Enlarged
Enhancement Virtual Element Methods (E\textsuperscript{2}VEM), designed to allow
the definition of a coercive bilinear form that involves only polynomial
projections. In this framework, it is not required to add an arbitrary
stabilizing bilinear form accounting for the non polynomial part of VEM
functions.  The method is based on the use of higher order polynomial
projections in the discrete bilinear form with respect to the standard one
\cite{Beirao2015b} and on a modification of the VEM space to allow the
computation of such projections. In particular, we extend the enhancement
property that is used in the definition of the VEM space (\cite{Ahmad2013},
\cite{Beirao2015b}). Indeed, the name of the method comes from this enlarged
enhancement property. The degree of polynomial enrichment is chosen locally on
each polygon, such that the discrete bilinear form is coercive,
and depends on the number of vertices of the polygon.  The resulting discrete
functional space has the same set of degrees of freedom of the one defined in
\cite{Beirao2015b}.

The proof of well-posedness is quite elaborate, thus in this paper we choose to
deal only with the lowest order formulation and, for the sake of simplicity, we
focus on the two dimensional Poisson's problem with homogenous Dirichlet
boundary conditions, the extension to general boundary conditions being
analogous to what is done for classical VEM. Moreover, the formulation and
proofs presented in this work can also be easily extended to the case of a non
constant anisotropic diffusion tensor. The extension to \fix{a} higher order
formulation will be the focus of an upcoming work\fix{, while in
  \cite{BBME2VEMLetter} a comparison between the proposed method and standard
  Virtual Elements from \cite{Beirao2015b} has been done, showing that the new
  formulation can speed up convergence in the case of anisotropic diffusion
  tensors. Indeed, since the stabilization term is isotropic by nature, its
  presence can induce larger errors.}

The outline of the paper is as follows. In \cref{sec:modelProb} we state our
model problem.  In \cref{sec:discrForm} we introduce the approximation
functional spaces and projection operators and we state the discrete
problem. \Cref{sec:wellposedness} contains the discussion about the
well-posedness of the discrete problem under suitable sufficient conditions on
the local projections. In \cref{sec:apriori} we prove optimal order
$\mathrm{H}^1$ a priori error estimates and in the supplementary materials the
$\lebl{}$ case. \fix{In \cref{sec:diffreact} we briefly discuss the extension of
  our approach to a diffusion-reaction model problem.} Finally,
\cref{sec:numericalresults} contains some numerical results assessing the rates
of convergence of the method.

Throughout the work, $(\cdot,\cdot)_{\omega}$ denotes the standard
$\mathrm{L}^2$ scalar product defined on a generic $\omega\subset\mathbb{R}^2$,
 $\trace{}{\partial \omega}{}$ denotes the trace operator, that restricts on the
boundary $\partial \omega$ an element of a space defined over
$\omega\subset\mathbb{R}^2$.  Inside the proofs, we decide to use a single
character $C$ for constants, independent of the mesh size, that appear in the
inequalities, which means that we suppose to take at each step the maximum of
the constants involved.

\section{Model Problem}
\label{sec:modelProb}
Let $\Omega \subset \mathbb{R}^2$ be a bounded open set. We are interested in solving the following
problem:
\begin{equation}
  \label{eq:modelProblem}
  \begin{cases}
    -\Delta U = f & \text{in $\Omega$,}
    \\
    U = 0 & \text{on $\partial \Omega$.}
  \end{cases}
\end{equation}
Defining $\a{}{}\colon \sobh[0]{1}{\Omega}\times \sobh[0]{1}{\Omega} \to \mathbb{R}$ such that,
\begin{equation}
  \label{eq:defBilFormA}
  \a{U}{W} := \scal[\Omega]{\nabla U}{\nabla W} \quad \forall U,W \in \sobh[0]{1}{\Omega},
\end{equation}
then, given $f\in\lebl{\Omega}$, the variational formulation of \cref{eq:modelProblem} is given by:
find $U\in \sobh[0]{1}{\Omega}$ such that,
\begin{equation}
  \label{eq:contVarForm}
  \a{U}{W} = \scal[\Omega]{f}{W} \quad \forall W\in\sobh[0]{1}{\Omega}\,.
\end{equation}



\section{Discrete formulation}
\label{sec:discrForm}

In order to define the discrete form of \cref{eq:contVarForm},
 $\Mh$ denotes a conforming polygonal tessellation of $\Omega$ and $E$ denotes a generic
polygon of $\Mh$. $\#\Mh$ denotes the number of polygons of $\Mh$ and the maximum diameter of all the polygons in $\Mh$ is denoted by $h$. Let
$\{x_i\}_{i=1}^{\NVE}$ be the $\NVE$ vertices of $E$, $\E{E}$ the set of its
edges and $\boldsymbol{n}^e=(n^e_x,n^e_y)$ the outward-pointing unit normal
vector to the edge $e$ of $E$.  We assume that $\Mh$ satifies the standard mesh
assumptions for VEM (see for instance
\cite{Beirao2017,Brenner2018}), i.e. $\exists \kappa > 0$ such that
\begin{enumerate}
\item for all $E\in \Mh$, $E$ is star-shaped with respect to a ball of radius $\rho\geq \kappa h_E$,
  where $h_E$ is the diameter of $E$;
\item for all edges $e\subset\partial E$, $\abs{e} \geq \kappa h_E$.
\end{enumerate}
Notice that the above conditions imply that, denoting by $\NVE$ the number of vertices of $E$, it
holds
\begin{equation}
  \label{eq:numVerticesBounded}
  \exists \NVmax>0\colon \forall E\in\Mh,\, \NVE \leq \NVmax \,.
\end{equation}

For any given $E\in\Mh$, let $\Poly{k}{E}$ be the space of polynomials of degree $k$ defined
on $E$. Let $\proj[\nabla]{1}{E}:\sobh{1}{E}\rightarrow\Poly{1}{E}$ be the $\sobh{1}{E}$-orthogonal
operator, defined up to a constant by the orthogonality condition: $\forall \,u\in\sobh{1}{E} $,
\begin{equation}\label{eq:PiNablaorthogonalitycondition}
  \scal[E]{\nabla\left(\proj[\nabla]{1}{E} u -u \right)}{\nabla p}=0 \;\; \forall \, p\in\Poly{1}{E}.
\end{equation}
In order to define $\proj[\nabla]{1}{E}$ uniquely, we choose any continuous and linear projection operator
$\mathrm{P}_0:\sobh{1}{E}\rightarrow\Poly{0}{E}$, whose continuity constant in $\sobh{1}{}$-norm is independent of $h_E$ and continuous with respect to deformations of the geometry, and we impose $\forall \,u\in\sobh{1}{E} $,
\begin{equation}\label{P0 definition}
  \mathrm{P}_0(\proj[\nabla]{1}{E} u -u) = 0.
\end{equation}
\begin{remark}
Under the current mesh assumptions, a suitable choice for $\mathrm{P}_0$ is the integral mean on the boundary of $E$, i.e.
\begin{equation*}
\mathrm{P}_0(u):=\frac{1}{\abs{\partial E}} \int_{\partial E}\limits \trace{u}{\partial E} \, ds \quad \forall u\in \sobh{1}{E}.
\end{equation*}
Notice that this is a common choice, see for instance \cite{Beirao2015b}.
\end{remark}

For any given $E\in\Mh$, let $l\in\mathbb{N}$ be given, as detailed in the next
section\fix{, where we will choose $l$ depending on $\NVE$ (see
  \cref{PiZeroEnhInjectivity})}. Let $\EN[E]{1,l}$ be the set of functions
$v\in\sobh{1}{E}$ satisfying
\begin{equation}
  \label{enhancement definition}
  \scal[E]{v}{p} =\scal[E]{\proj[\nabla]{1}{E} v}{p} \; \forall p\in\Poly{l+1}{E} \,.
\end{equation}
We define the Enlarged Enhancement Virtual Space of order $1$ as
\begin{equation*}
  \V[E]{1,l}:=\{ v\in \EN[E]{1,l}: \Delta v\in \Poly{l+1}{E}, \;\; \trace{v}{e}{}\in \Poly{1}{e}\;\forall e \in\E{E}, \;\; v\in\cont{\partial E}\}\,.
\end{equation*}
We define as degrees of freedom of this space the values of functions at the vertices of $E$ (see \cite{Beirao2013a,Beirao2015b}).

Moreover, let $\boldell\in \mathbb{N}^{\#\Mh}$ be a vector
 and $\boldell(E)$ denote the element corresponding
to the polygon $E$, we define the global discrete space as
\begin{equation*}
  \V{1,\boldell} := \{v\in\sobh[0]{1}{\Omega}\colon v_{|E} \in \V[E]{1,l},\,\text{where $l = \boldell(E)$}\}
  \,.
\end{equation*}
Note that $v\in\V{1,\boldell}$ is a continuous function that is a polynomial of degree $1$ on each edge of the mesh.

To define our discrete bilinear form, let
$\proj{l}{E} \nabla: \sobh{1}{E} \rightarrow \PolyDouble{l}{E} $ be the $\lebl{E}$-projection operator of the
gradient of functions in $\sobh{1}{E}$, defined, $\forall u\in \sobh{1}{E}$, by the orthogonality condition
\begin{equation}\label{eq:PiZeroLGrad orthogonality condition}
  \scal[E]{\proj{l}{E} \nabla u}{\boldsymbol{p}} = \scal[E]{\nabla u}{\boldsymbol{p}} \;\; \forall \boldsymbol{p} \in \left[\Poly{l}{E}\right]^2 \,.
\end{equation}
\begin{remark}
For each function $u\in \V[E]{1,l}$, the above projection is computable given the degrees of freedom of $u$, applying the Gauss-Green formula and exploiting \cref{enhancement definition}.
\end{remark}

Let $\ahE{}{}\colon \V[E]{1,l}\times \V[E]{1,l} \to \mathbb{R}$ be defined as
\begin{equation}
  \label{eq:defahE}
  \ahE{u}{v} := \scal[E]{\proj{l}{E} \nabla u}{\proj{l}{E} \nabla v} \quad \forall u,v \in \V[E]{1,l}\,,
\end{equation}
and $\a[h]{}{}\colon \V{1,\boldell}\times \V{1,\boldell} \to \mathbb{R}$ as
\begin{equation}
  \label{eq:defah}
  \a[h]{u}{v} := \sum_{E\in\Mh}\ahE{u}{v} \quad \forall u,v \in \V{1,\boldell} \,.
\end{equation}
We can state the discrete problem as: find $u\in\V{1,\boldell}$ such that
\begin{equation}
  \label{eq:discrVarForm}
  \ah{u}{v} = \sum_{E\in\Mh}\scal[E]{f}{\proj{0}{E}v}
  \quad \forall v \in \V{1,\boldell}\,,
\end{equation}
where, $\forall E\in\Mh$, $\proj{0}{E}\colon \lebl{E} \to \mathbb{R}$ is the $\lebl{E}$-projection, defined by
\begin{equation}
  \label{eq:defIntegralMean}
  \proj{0}{E} v := \frac{1}{\abs{E}} \scal[E]{v}{1} \quad \forall v\in\ourfix{\lebl{E}} \,.
\end{equation}
The above projection is computable for any given $v\in\V[E]{1,l}$ exploiting \cref{enhancement definition}.


\section{Well-posedness}
\label{sec:wellposedness}

This section is devoted to prove the well-posedness of the discrete problem
stated by \cref{eq:discrVarForm}, under suitable sufficient conditions on
$\boldell$. The main result is given by \cref{PiZeroEnhInjectivity}, that
induces the existence of an equivalent norm on $\V{1,\boldell}$, which implies
the well-posedness of \cref{eq:discrVarForm}.

First, we define, for any given $l\in\mathbb{N}$,
\begin{equation}
  \label{eq:defBadPolynomials}
  \BadPoly{l}{E} = \left\{\boldsymbol{p}\in\PolyDouble{l}{E}\colon
    \int_{\partial E}\boldsymbol{p}\cdot \boldsymbol{n}^{\partial E}
    \trace{v-\mathrm{P}_0(v)}{\partial E}{} = 0
    \quad \forall v \in \V[E]{1,l}\right\} \,.
\end{equation}

\ourfixBis{%
  Notice that the dimension of $\BadPoly{l}{E}$ generally depends on the
  geometry of the polygon and the definition of $\mathrm{P}_0$, but in
  \cref{th:badpolyupperbound} we provide an upper bound for
  $\dim\BadPoly{l}{E}$. }

Then, the following result holds.

\begin{theorem} \label{PiZeroEnhInjectivity} Let $E\in\Mh$, $u\in \V[E]{1,l}$
  and $l\in\mathbb{N}$ such that \fix{the following condition is satisfied:}
  \begin{gather}
    \label{eq:suffCondWellPos}
    \fix{(l+1)(l+2) - \dim \BadPoly{l}{E} \geq \NVE - 1},
  \end{gather}
  then
  \begin{equation} \label{relation of injectivity to prove} \proj{l}{E}\nabla
    u=0 \imply {\nabla u}_{|_E}=0.
  \end{equation}
\end{theorem}

We omit in the following the proof of the case of triangles ($\NVE=3$ and
$l=0$), indeed \fix{if $E$ is a triangle, $\V[E]{1,l} = \Poly{1}{E}$
  $\forall l\geq 0 $, and then $\proj{l}{E}\nabla u = \nabla u$
  $\forall l\geq 0 $. Moreover, an explicit computation yields
  $\dim\BadPoly{0}{E}=0$ if $E$ is a triangle}. Then, for technical reasons, the
proof of \cref{PiZeroEnhInjectivity} \fix{for a general polygon} is split into
two results, described in \cref{sec:equivalent-inf-sup} and in
\cref{sec:fortintrick}, respectively. The proof relies on an auxiliary inf-sup
condition that is proved by constructing a suitable Fortin operator, whose
existence is guaranteed under condition \cref{eq:suffCondWellPos}.


\subsection{Auxiliary inf-sup condition}
\label{sec:equivalent-inf-sup}

In this section, after some auxiliary results, we prove through 
\cref{prop: Equivalence with inf-sup} that \cref{relation of injectivity to
  prove} is satisfied if the auxiliary inf-sup condition \cref{inf-sup
  condition on b} holds true.

\begin{lemma} \label{PiNabla constant} Let $u\in \V[E]{1,l}$, with $l\geq 0$.
  Then
  \begin{equation*}
    \proj{l}{E} \nabla u =0 \imply \proj[\nabla]{1}{E} u \in \Poly{0}{E}.
  \end{equation*}
\end{lemma}
\begin{proof}
  Applying \cref{eq:PiZeroLGrad orthogonality condition}, we have
  \begin{equation*}
    \proj{l}{E} \nabla u =0 \imply \scal[E]{\nabla u}{\boldsymbol{p}} = 0 \;\; \forall \boldsymbol{p} \in \left[\Poly{l}{E}\right]^2,
  \end{equation*}
  that implies
  \begin{equation} \label{eq:gradugradp} \scal[E]{\nabla u}{\nabla p} = 0 \;\;
    \forall p \in \Poly{1}{E},
  \end{equation}
  thanks to the relation
  $\nabla\Poly{1}{E}\subseteq\nabla\Poly{l+1}{E}\subseteq\left[\Poly{l}{E}\right]^2$.
  Given \cref{eq:gradugradp} and \cref{eq:PiNablaorthogonalitycondition},
  \begin{align*}
    \scal[E]{\nabla\proj[\nabla]{1}{E} u}{\nabla p} = 0 \;\; \forall p \in \Poly{1}{E} &\imply \nabla\proj[\nabla]{1}{E} u = 0 \\ & \imply \proj[\nabla]{1}{E} u \in \Poly{0}{E}.
  \end{align*}
\end{proof}
\begin{lemma}
  Let $u\in \V[E]{1,l}$. If $\proj{l}{E} \nabla u =0$, then \cref{enhancement
    definition} can be rewritten as
  \begin{equation}\label{eq:enhancement rewrittening}
    \scal[E]{u}{p} =\mathrm{P}_0\left(u\right)\cdot\scal[E]{1}{p} \; \forall p\in\Poly{l+1}{E},
  \end{equation}
  where $\mathrm{P}_0$ is the projection operator chosen in
  \cref{sec:discrForm}.
\end{lemma}
\begin{proof}
  Applying  \cref{PiNabla constant} and \cref{P0 definition},
  \begin{equation*}
    \proj{l}{E} \nabla u =0 \imply \proj[\nabla]{1}{E} u = \mathrm{P}_0\left(u\right).
  \end{equation*}
  Then, \cref{enhancement definition} provides \cref{eq:enhancement
    rewrittening}.
\end{proof}

We now need to introduce some \fix{notation} and definitions. First, let $\Triang{E}$ denote the sub-triangulation of $E$ obtained linking each
vertex of $E$ to the centre of the ball with respect to which $E$ is star-shaped, denoted by $x_C$.  Let us define the set of internal edges of
the triangulation $\Triang{E}$ as $\intEdges$.  For any $i=1,\ldots, \NVE$, let $\tau_i\in \Triang{E}$ be
the triangle whose vertices are $x_i$, $x_{i+1}$ and $x_C$.  Let $e_i$ denote
the edge $\overrightarrow{x_Cx_i}\in\intEdges$ and by $\boldsymbol{n}^{e_i}$ the
outward-pointing unit normal vector to the edge $e_i$ of $\tau_i$.
 Then, for each polygon $E$, we can define the reference polygon $\hat{E}$,
 such that the mapping $F : \hat{E}\rightarrow E $ is given by
 \begin{equation}
   \label{def:affinemapping}
   x=h_E\hat{x}+x_C.
 \end{equation}
Let $\Sigma$ be the set of all admissible reference polygons, i.e. satisfying the mesh assumptions with the same regularity parameter as the polygons in the mesh.
 \begin{lemma}[{\cite[Proof of Lemma 4.9]{Cangiani2016}}]
   \label{lem:compactnessSigma}
   $\Sigma$ is compact.
 \end{lemma}
 
\begin{definition}
  Let $\HT{E}$ be the broken Sobolev space
  \begin{equation*}
    \HT{E}:= \left\{v\colon \eval{v}{\tau}\in\sobh{1}{\tau} \; \forall \tau\in\Triang{E}\right\}.
  \end{equation*}
  Let $u\in\HT{E}$, we define $\forall e_i\in\intEdges$ the jump function
  $\jmp[e_i]{\cdot}:\HT{E}\rightarrow \lebl{e_i}$ such that
  \begin{equation*}
    \jmp[e_i]{u}:= \trace{u_{|_{\tau_i}}}{e_i}{} - \trace{u_{|_{\tau_{i-1}}}}{e_i}{}.
  \end{equation*}
  Moreover, $\jmp[\intEdges]{u}$ denotes the vector containing the jumps of $u$
  on each $e_i\in\intEdges$.  We endow $\HT{E}$ with the following seminorm and
  norm : $\forall u \in \HT{E},$
  \begin{align}
    \quad &\seminorm[\HT{E}]{u}^2
            :=\sum_{\tau\in\Triang{E}}\norm[\lebldouble{\tau}]{\nabla u}^2 +\sum_{i=1}^{\NVE}\limits \norm[\lebl{e_i}]{\jmp[e_i]{u}}^2,
    \\
    \label{broken norm of H1}
          &\norm[\HT{E}]{u} ^2 :=\seminorm[\HT{E}]{u}^2+
            \sum_{\tau\in\Triang{E}}\norm[\lebl{\tau}]{ u}^2 \,.
  \end{align}
\end{definition}
\ourfixBis{\begin{definition}
Let $\spaceV$ be given by
  \begin{equation}\label{definition of V space}
    \spaceV:=\{\bs{v}\in\Hdiv{\tau} \;\forall \tau\in\Triang{E} :  \forall e_i\in\intEdges,\,\jmp[e_i]{\bs{v}}\in\lebl[\infty]{e_i}\}.
  \end{equation}
  Then $\forall \bs v\in \spaceV$, we define its seminorm and its norm:
  \begin{align*}
    &\seminorm[\spaceV]{\bs v}^2:=\sum_{\tau\in\Triang{E}} \norm[\lebl{\tau}]{\div{\bs{v}}}^2 + h_E^2\norm[\leblinfI]{\jmp[\intEdges]{\bs{v}}}^2,\\
    &\norm[\spaceV]{\bs v}^2:= \seminorm[\spaceV]{\bs v}^2 + \sum_{\tau\in\Triang{E}} \norm[\lebldouble{\tau}]{\bs v}^2
  \end{align*}
  where
  \begin{equation*}
    \norm[\leblinfI]{\jmp[\intEdges]{\bs v}} := \max_{i=1,\ldots,\NVE} \limits \norm[{\lebl[\infty]{e_i}}]{\jmp[e_i]{\bs v}}.
  \end{equation*}
\end{definition}}
\begin{remark}
  Let us observe that $\ourfixBis{\PolyDouble{l}{E}} \subset \spaceV$.
  Hence, we can use $\norm[\ourfixBis{\spaceV}]{\cdot}$ as a norm for $\PolyDouble{l}{E}$.
  Notice that, since $\PolyDouble{l}{E}\subset[\cont{E}]^2$,
  $\norm[\leblinfI]{\jmp[\intEdges]{\boldsymbol{p}}} =
  0,\;\forall\boldsymbol{p}\in \PolyDouble{l}{E}$ .
\end{remark}

\begin{definition}
  Let $\V[E,\mathrm{P}_0]{1,l}$ be the space
  \begin{equation}
    \V[E,\mathrm{P}_0]{1,l} := \left\{ v\in \V[E]{1,l}: \mathrm{P}_0(v)=0 \right\}.
  \end{equation}
\end{definition}

\begin{definition}
  Denoting by $\{\psi_i\}_{i=1}^{\NVE}$ the set of Lagrangian basis functions of $\V[E]{1,l}$, let
  $\mathcal{Q} (\partial E)$ be the vector space
  \begin{equation}\label{definition of Q}
    \mathcal{Q} (\partial E) := \mathrm{span} \left\{ \trace{\psi_i-\mathrm{P}_0(\psi_i)}{\partial E}{}\right\}, \quad \forall i = 1,\ldots,\NVE -1.
  \end{equation}
  \fix{We remark that the above space is made up of continuous piecewise linear
    polynomials on each edge.}  Notice that
  $\forall q\in \mathcal{Q} (\partial E)$,
  $\exists !v\in\V[E,\mathrm{P}_0]{1,l}$ such that $q=\trace{v}{\partial E}{}$.
\end{definition}

\begin{definition}
  Let $\RQE$ be the vector space, lifting of $\mathcal{Q}(\partial E)$ on $E$,
  given by:
  \begin{equation} \label{Definition of RQE} \RQE:=\left\{
      \bar{q}_{|\tau}\in\Poly{1}{\tau}\, \forall \tau\in\Triang{E},\,
      \trace{\bar{q}}{\partial E}{} \in \mathcal{Q}(\partial
      E),\,\bar{q}(x_C)=0\right\}.
  \end{equation}
  We note that $\RQE\subset \HT{E}\cap\cont{E}$.  Hence, we use the norm
  $\norm[\HT{E}]{\cdot}$ defined in \cref{broken norm of H1} as a norm for
  $\RQE$.  Notice that
  $\sum_{i=1}^{\NVE}\limits \norm[\lebl{e_i}]{\jmp[e_i]{\bar{q}}} = 0 $. Let
  $\{r_j\}_{j=1}^{\NVE-1}$ denote a basis of $\RQE$.
\end{definition}

Now, we can introduce the bilinear form $b$ which is used in 
\cref{prop: Equivalence with inf-sup}.
\begin{definition}
  Let $b:\RQE\times\ourfixBis{\spaceV}\rightarrow\mathbb{R}$, such that
  $\forall \bar{q}\in\RQE,\,\forall\,\boldsymbol{v}\in\ourfixBis{\spaceV}$
  \begin{equation}
    \label{eq:defB}
    b(\bar{q},\boldsymbol v):=\int_{\partial E} \limits \bar{q}\,\boldsymbol{v}\cdot n^{\partial E} \, ds.
  \end{equation}
  Applying the divergence theorem, we can rewrite the form $b$:
  \begin{equation}
    \label{eq:BAfterDivTheorem}
    b(\bar{q},\boldsymbol v) = \sum_{\tau\in\Triang{E}} \limits \int_{ \tau} \left[\nabla\bar{q}\,\boldsymbol{v} + \bar{q}\,\div \boldsymbol{v} \right] \, dx - \, \sum_{i=1}^{\NVE} \int_{e_i} \limits \trace{\bar{q}}{e_i}{}\jmp[e_i]{\boldsymbol{v}}\cdot \boldsymbol{n}^{e_i} ds.
  \end{equation}
\end{definition}


The following lemma gives the continuity of the bilinear form $b$.
\begin{lemma}\label{lem: continuity of b}
  Let $b$ be given by \cref{eq:defB}. Then $b$ is a bilinear form and
  $\exists \fix{C_b}>0$ independent of $h_E$ such that
  \begin{equation*}
    b(\bar{q},\boldsymbol v)\leq \fix{C_b}\norm[\HT{E}]{\bar{q}}\norm[\ourfixBis{\spaceV}]{\boldsymbol{v}} \;\; \forall \bar{q}\in\RQE,\,\forall\,\boldsymbol{v}\in\ourfixBis{\spaceV} \,.
  \end{equation*}
\end{lemma}
\begin{proof}
The proof of this lemma can be found in the supplementary materials of this paper.
\end{proof}

The following proposition is the first step towards the proof of \cref{PiZeroEnhInjectivity}.
\begin{proposition} \label{prop: Equivalence with inf-sup}
Let $b$ the continuous bilinear form defined by \cref{eq:defB}.
If $\exists \beta>0 $, independent of $h_E$, such that
\begin{equation} \label{inf-sup condition on b}
\forall \bar{q}\in \RQE, \quad \sup_{\boldsymbol{p} \in \PolyDouble{l}{E}} \limits\frac{b(\bar{q},\boldsymbol{p})}{\norm[\ourfixBis{\spaceV}]{\boldsymbol{p}}} \geq \beta\norm[\HT{E}]{\bar{q}} \,,
\end{equation}
then \cref{relation of injectivity to prove} holds true.
\end{proposition}
\begin{proof}
Let $u\in\V[E]{1,l}$, given \cref{eq:PiZeroLGrad orthogonality condition},
\begin{equation*}
\proj{l}{E}\nabla u=0\imply\scal[E]{\nabla u}{\boldsymbol{p}} = 0 \;\; \forall \boldsymbol{p} \in \left[\Poly{l}{E}\right]^2.
\end{equation*}
Applying Gauss-Green formula, the previous relation becomes
\begin{equation*}
\scal[E]{\nabla u}{\boldsymbol{p}} = \scal[\partial E]{\trace{u}{\partial E}{}}{\boldsymbol{p}\cdot n^{\partial E}} - \scal[E]{u}{\div \boldsymbol{p}} = 0 \;\; \forall \boldsymbol{p}\in \left[\Poly{l}{E}\right]^2.
\end{equation*}
Since $\div \boldsymbol{p}\in\Poly{l-1}{E}$ we apply \cref{eq:enhancement rewrittening} and we obtain
\begin{equation*}
\scal[\partial E]{\trace{u}{\partial E}{}}{\boldsymbol{p}\cdot n^{\partial E}} - \mathrm{P}_0\left(u\right)\cdot\scal[E]{1}{\div \boldsymbol{p}} = 0 \;\; \forall \boldsymbol{p} \in \left[\Poly{l}{E}\right]^2.
\end{equation*}
Then we can apply the divergence theorem and find the relation
\begin{equation} \label{condition restricted on the boundary}
\scal[\partial E]{\trace{u-\mathrm{P}_0\left(u\right)}{\partial E}{}}{\boldsymbol{p}\cdot n^{\partial E}} = 0 \;\; \forall \boldsymbol{p} \in \left[\Poly{l}{E}\right]^2.
\end{equation}
We have $q=\trace{u-\mathrm{P}_0\left(u\right)}{\partial E}{}\in\mathcal{Q}(\partial E)$ ($\mathcal{Q}(\partial E)$  defined in \cref{definition of Q}).
Let $\bar{q}\in\RQE$ be the lifting of $q$ ($\RQE$ defined in \cref{Definition of RQE}), then the relation \cref{condition restricted on the boundary} is
\begin{equation*}
b(\bar{q},\boldsymbol{p}) =0 \;\; \forall \boldsymbol{p} \in \left[\Poly{l}{E}\right]^2.
\end{equation*}
Then, since $b$ is a continuous bilinear form, \cref{inf-sup condition on b} implies $q\equiv 0$.
Finally, since $u\in\V[E]{1,l}$, then $u=\mathrm{P}_0\left(u\right)$.
\end{proof}

\subsection{Proof of the inf-sup condition}
\label{sec:fortintrick}

In this section we show that \cref{inf-sup condition on b} holds with $\beta$
independent of $h_E$. The proof relies on the technique known as Fortin trick
\cite{Brezzi_Mixed}, that consists in the following classical result.

\begin{proposition}[{\cite[Proposition 5.4.2]{Brezzi_Mixed}}] \label{prop:
    Brezzi _1} Assume that there exists an operator
  $\Pi_E:\ourfixBis{\spaceV}\rightarrow\PolyDouble{l}{E}$ that satisfies,
  $\forall \boldsymbol{v} \in \ourfixBis{\spaceV}$,
  \begin{equation}
    \label{eq:b-compatibility}
    b(\bar{q},\Pi_E\boldsymbol{v}-\boldsymbol{v})= 0 \;\; \forall \bar{q}\in\RQE \,,
  \end{equation}
  and assume that there exists a constant $C_\Pi>0$, independent of $h_E$, such
  that
  \begin{equation}
    \label{eq:Fortincontinuity}
    \norm[\ourfixBis{\spaceV}]{\Pi_E\boldsymbol{v}}\leq C_\Pi\norm[\ourfixBis{\spaceV}]{\boldsymbol{v}} \;\; \forall \boldsymbol{v}\in \ourfixBis{\spaceV}.
  \end{equation}
  Assume moreover that $\exists \eta > 0$, independent of $h_E$ such that
  \begin{equation}
    \label{eq:infsup-cont}
    \inf_{q \in \RQE} \sup_{\boldsymbol{v} \in \ourfixBis{\spaceV}}
    \frac{b(q, \boldsymbol{v})}{
      \norm[\HT{E}]{q}
      \norm[\ourfixBis{\spaceV}]{\boldsymbol{v}}
    } \geq \eta \,.
  \end{equation}
  Then the discrete inf-sup condition \cref{inf-sup condition on b} is
  satisfied, with $\beta = \frac{\eta}{C_\Pi}$.
\end{proposition}

\begin{remark}
  The inf-sup constant $\beta$ in \cref{inf-sup condition on b} has to be
  independent of the mesh size in order to guarantee that the constant in
  \cref{eq:lowerBoundBoldellNorm}, involved in the coercivity of the bilinear
  form of \cref{eq:discrVarForm}, is independent of the mesh size.
\end{remark}
\begin{remark}
  The operator $\Pi_E$ defined in the following is such that the constant
  $C_\Pi$ depends on $\NVmax$ and on the continuity constant of $\mathrm{P}_0$,
  both are bounded independently of $h_E$ by assumption.
\end{remark}

Following the above results, we have to prove \cref{eq:infsup-cont} and to show
the existence of the operator $\Pi_E$ satisfying \cref{eq:b-compatibility} and
\cref{eq:Fortincontinuity}. In the following proposition we achieve the first
task.

\begin{proposition}
  \label{prop:infsup-cont}
  Let $b\colon \RQE\times\ourfixBis{\spaceV}\rightarrow\mathbb{R}$ be defined by
  \cref{eq:defB}. Then the inf-sup condition \cref{eq:infsup-cont} holds true.
\end{proposition}
\begin{proof}
  \ourfixBis{ Let $q\in\RQE$ be given arbitrarily.  For any
    $\tau_i\in\Triang{E}$, we recall that the vertices of $\tau_i$ are $x_C$,
    $x_i=\left(\begin{smallmatrix}x_{i,1}\\x_{i,2}\end{smallmatrix}\right)$ and
    $x_{i+1}=\left(\begin{smallmatrix}x_{i+1,1}\\x_{i+1,2}\end{smallmatrix}\right)$.
    Let $e^{\partial}_i$ and $\bs{n}^{e^{\partial}_i}$ denote the edge
    $\overrightarrow{x_ix_{i+1}}$ and the outward-pointing unit normal vector to
    the edge $e^{\partial}_i$, respectively.  Let $\bs\varphi^1_i,$
    $\bs\varphi^2_i\in\PolyDouble{1}{\tau_i}$ be given such that
    \begin{equation}
      \bs\varphi^1_i(x_1,x_2)= \begin{pmatrix} x_2-x_{i,\,2}\\-(x_1-x_{i,\,1}) \end{pmatrix}\,, 
      \bs\varphi^2_i(x_1,x_2) =\abs{e^{\partial}_i} \bs{n}^{e^{\partial}_i} =
      \begin{pmatrix}
        x_{i+1,2}-x_{i,2} \\ -(x_{i+1,1}-x_{i,1})
      \end{pmatrix}
      \,.
    \end{equation}
    Let
    $\mathrm{B}^{\partial}(\tau_i):=\mathrm{span}\{\bs\varphi^1_i,\bs\varphi^2_i\}\subset\PolyDouble{1}{\tau_i}$.
    Notice that $\forall \bs{v}\in\mathrm{B}^{\partial}(\tau_i)$ $\div\bs v=0$
    and $\norm[\lebl{e^{\partial}_i}]{\bs v\cdot \bs n^{e^{\partial}_i}}$ is a
    norm on $\mathrm{B}^{\partial}(\tau_i)$. Indeed, if
    $\bs v \in \mathrm{B}^{\partial}(\tau_i)$ and
    $\norm[\lebl{e^{\partial}_i}]{\bs v\cdot \bs n^{e^{\partial}_i}}=0$, then an
    explicit computation yields $\bs v \equiv \bs{0}$ on $\tau_i$.  Moreover,
    let
    $\mathrm{B}^{\partial}(\Triang{E}):=\{\bs v:\bs
    v_{|_{\tau}}\in\mathrm{B}^{\partial}(\tau) \; \forall \tau \in
    \Triang{E}\}\subset\spaceV$.  Notice that
    $\forall \bs v \in \mathrm{B}^{\partial}(\Triang{E})$
    $\sum_{\tau \in\Triang{E}}\norm[\lebl{\tau}]{\div \bs v}^2 = 0$ and
    $\norm[\lebl{\partial E}]{\bs v\cdot \bs n^{\partial E}}$ is a norm on
    $\mathrm{B}^{\partial}(\Triang{E})$.  We define
    $\bs v^\star\in\mathrm{B}^{\partial}(\Triang{E})$ such that
    $\eval{\boldsymbol{v}^\star}{ \partial E} \cdot \boldsymbol{n}^{\partial E}
    = \trace{q}{\partial E}{}$ .  In particular, $\forall \tau_i\in\Triang{E}$
    $\eval{\boldsymbol{v}^\star}{\tau_i}
    =\frac{q(x_{i+1})-q(x_i)}{\abs{e^{\partial}_i}}\bs\varphi^1_i+\frac{q(x_i)}{\abs{e^{\partial}_i}}\bs\varphi^2_i$.}
  Notice that
  $\norm[\lebl{\partial E}]{q} = \norm[\lebl{\partial
    E}]{\boldsymbol{v}^\star\cdot \boldsymbol{n}^{\partial E}}$.Then,
  \begin{equation}
    \label{eq:infsup-cont_firststep}
    \sup_{\boldsymbol{v}\in\ourfixBis{\spaceV}}\frac{b(q, \boldsymbol{v})}{
      \norm[\HT{E}]{q}
      \norm[\ourfixBis{\spaceV}]{\boldsymbol{v}}} \geq
    \frac{b(q, \boldsymbol{v}^\star)}{
      \norm[\HT{E}]{q}
      \norm[\ourfixBis{\spaceV}]{\boldsymbol{v}^\star}
    } = \frac{\norm[\lebl{\partial E}]{q}}{\norm[\HT{E}]{q}}
    \frac{ \norm[\lebl{\partial E}]{\boldsymbol{v}^\star\cdot \boldsymbol{n}^{\partial E}}}{
      \norm[\ourfixBis{\spaceV}]{\boldsymbol{v}^\star}
    }  \,.
  \end{equation}
  We have to estimate from below the last two factors. We notice that
  $\norm[\lebl{\partial E}]{q}$ is a norm on $\RQE$\fixBis{, since
    $q \in\Poly{1}{\Triang{E}}$ and $q(x_C) = 0$}. Thus, we can exploit the
  equivalence of norms on finite dimensional spaces. Hence, regarding the first
  norm, we get, by a scaling argument,
  \begin{equation}
    \label{eq:qnormestim}
    \begin{split}
      \norm[\lebl{\partial E}]{q}^2
      &= \sum_{e\in\partial E}
        \norm[\lebl{e}]{q}^2 = h_E\sum_{\hat{e}\in\partial \hat{E}}
        \norm[\lebl{\hat{e}}]{\hat{q}}^2 \geq Ch_E
        \left(\sum_{\hat{\tau}\in\Triang{\hat{E}}}\limits
        \norm[\lebl{\hat{\tau}}]{\hat{q}}^2 +
        \norm[\lebldouble{\hat{\tau}}]{\hat{\nabla} \hat{q}}^2\right)
      \\
      &= Ch_E \left(\sum_{\tau\in\Triang{E}} h_E^{-2}\norm[\lebl{\tau}]{q}^2 +
        \norm[\lebldouble{\tau}]{\nabla q}^2\right) \geq Ch_E
        \min\{1,h_E^{-2}\} \norm[\HT{E}]{q}^2
      \\
      & \geq C h_E \norm[\HT{E}]{q}^2\,.
    \end{split}
  \end{equation}
  Notice that the constant above is independent of the choice of reference
  element by \cref{lem:compactnessSigma}. The second norm is estimated using the
  definition of dual norm and the trace inequality
  \begin{equation*}
    \norm[\lebl{\partial E}]{\trace{w}{\partial E}{}} \leq
    Ch_E^{\frac12} \left(h_E^{-2}\norm[\lebl{E}]{w}^2 + \norm[\lebldouble{E}]{\nabla w}^2\right)^{\frac12}
    \quad \forall w \in\sobh{1}{E} \,,
  \end{equation*}
  as follows:
  \begin{equation}
    \label{eq:vstarnnormestim_1}
    \begin{split}
      \norm[\lebl{\partial E}]{\boldsymbol{v}^\star\cdot
      \boldsymbol{n}^{\partial E}}
      &= \sup_{\chi\in\lebl{\partial
        E}}\frac{\scal[\partial E]{\boldsymbol{v}^\star\cdot
        \boldsymbol{n}^{\partial E}}{\chi}}{\norm[\lebl{\partial
        E}]{\chi}}\geq %
        \sup_{w\in\sobh{1}{E}}\frac{\scal[\partial E]{\boldsymbol{v}^\star\cdot
        \boldsymbol{n}^{\partial E}}{\trace{w}{\partial
        E}{}}}{\norm[\lebl{\partial E}]{\trace{w}{\partial E}{}}}
      \\
      &\geq C h_E^{-\frac12} \ourfixBis{h_E} \sup_{w\in\sobh{1}{E}}
        \frac{\ourfixBis{ h_E^{-1}\scal[\partial E]{\boldsymbol{v}^\star\cdot \boldsymbol{n}^{\partial E}}
        {\trace{w}{\partial E}{}}}}{\left(h_E^{-2}\norm[\lebl{E}]{w}^2
        + \norm[\lebldouble{E}]{\nabla w}\right)^{\frac12}}
        \,.
    \end{split}
  \end{equation}
  Let $w^\star\in \sobh{1}{E}$ be such that
  \begin{equation*}
    \scal[E]{\nabla w^\star}{\nabla \varphi} + h_E^{-2} \scal[E]{w^\star}{\varphi} =
    \ourfixBis{ h_E^{-1}\scal[\partial E]{\boldsymbol{v}^\star\cdot
        \boldsymbol{n}^{\partial E}}{\trace{w}{\partial E}{}}} \,.
  \end{equation*}
  Notice that $\hat{w}^\star = w^\star \circ F$ ($F$ being the mapping defined
  by \cref{def:affinemapping}) is the solution of
  \begin{equation*}
    \scal[\hat{E}]{\hat{\nabla} \hat{w}^\star}{\hat{\nabla} \hat{\varphi}}
    + \scal[\hat{E}]{\hat{w}^\star}{\hat{\varphi}}
    =
    \scal[\partial \hat{E}]{\hat{\boldsymbol{v}}^\star\cdot \boldsymbol{n}^{\partial \hat{E}}}
    {\trace{\hat{\varphi}}{\partial \hat{E}}{}}
    \quad \forall \hat{\varphi} \in\sobh{1}{\hat{E}} \,.
  \end{equation*}
  Notice that
  \begin{equation*}
    \sup_{\hat{w}\in\sobh{1}{\hat{E}}}\frac{\scal[\partial
      \hat{E}]{\hat{\boldsymbol{v}}^\star\cdot \boldsymbol{n}^{\partial \hat{E}}}{\trace{\hat{w}}{\partial
          \hat{E}}{}}}{\norm[\sobh{1}{\hat{E}}]{\hat{w}}} = \frac{\scal[{\partial
        \hat{E}}]{\hat{\boldsymbol{v}}^\star\cdot \boldsymbol{n}^{\partial \hat{E}}}{\trace{\hat{w}^\star}{\partial
          \hat{E}}{}}}{\norm[\sobh{1}{\hat{E}}]{\hat{w}^\star}} \,.
  \end{equation*}
  This relation holds true since the \textit{greater than inequality} is trivial
  using the definition of $\sup$ and the \textit{less than inequality} can be
  proved applying the property of inner products
  $\abs{\scal{x}{y}}^2\leq\scal{x}{x}\scal{y}{y}$,
  indeed
  \begin{align*}
    \sup_{\hat{w}\in\sobh{1}{\hat{E}}}\frac{\scal[\partial
    \hat{E}]{\hat{\boldsymbol{v}}^\star\cdot \boldsymbol{n}^{\partial \hat{E}}}{\trace{\hat{w}}{\partial
    \hat{E}}{}}}{\norm[\sobh{1}{\hat{E}}]{\hat{w}}}
    &=
      \sup_{\hat{w}\in\sobh{1}{\hat{E}}}\frac{\scal[\hat{E}]{\hat{\nabla} \hat{w}^\star}{\hat{\nabla} \hat{w}}
      + \scal[\hat{E}]{\hat{w}^\star}{\hat{w}}}{\norm[\sobh{1}{\hat{E}}]{\hat{w}}}
    \\
    &\leq \sup_{\hat{w}\in\sobh{1}{\hat{E}}} \frac{\norm[\sobh{1}{\hat{E}}]{\hat{w}^\star}\norm[\sobh{1}{\hat{E}}]{\hat{w}}}{\norm[\sobh{1}{\hat{E}}]{\hat{w}}}
    \\ 
    &=\frac{\norm[\sobh{1}{\hat{E}}]{\hat{w}^\star}^2}{\norm[\sobh{1}{\hat{E}}]{\hat{w}^\star}} 
      = 
           \frac{\scal[{\partial
           \hat{E}}]{\hat{\boldsymbol{v}}^\star\cdot \boldsymbol{n}^{\partial \hat{E}}}{\trace{\hat{w}^\star}{\partial
           \hat{E}}{}}}{\norm[\sobh{1}{\hat{E}}]{\hat{w}^\star}} \,.
    \end{align*}
    Then, by choosing $\varphi=w^\star$ and $\hat{\varphi}=\hat{w}^\star$ in the
    equations above we get
    \begin{equation*}
      \begin{split}
        \sup_{w\in\sobh{1}{E}}\frac{\ourfixBis{ h_E^{-1}\scal[\partial E]{\boldsymbol{v}^\star\cdot
        \boldsymbol{n}^{\partial E}}{\trace{w}{\partial E}{}}}}{\left(h_E^{-2}\norm[\lebl{E}]{w}^2 + \norm[\lebldouble{E}]{\nabla w}\right)^{\frac12}}
        & \geq
          \frac{h_E^{-2}\norm[\lebl{E}]{w^\star}^2 + \norm[\lebldouble{E}]{\nabla
          w^\star}^2}{\left(h_E^{-2}\norm[\lebl{E}]{w^\star}^2 + \norm[\lebldouble{E}]{\nabla w^\star}\right)^{\frac12}}
        \\
        &= \left( h_E^{-2}\norm[\lebl{E}]{w^\star}^2 +
          \norm[\lebldouble{E}]{\nabla w^\star}^2 \right)^{\frac12}
        \\
        &=  \left( \norm[\lebl{\hat{E}}]{\hat{w}^\star}^2 +
          \norm[\lebldouble{\hat{E}}]{\hat{\nabla} \hat{w}^\star}^2 \right)^{\frac12} \\
        &=  \frac{\scal[{\partial \hat{E}}]{\hat{\boldsymbol{v}}^\star\cdot
          \boldsymbol{n}^{\partial \hat{E}}}{\trace{\hat{w}^\star}{\partial
          \hat{E}}{}}}{\norm[\sobh{1}{\hat{E}}]{\hat{w}^\star}}
        \\
        &=  \sup_{\hat{w}\in\sobh{1}{\hat{E}}}\frac{\scal[\partial
          \hat{E}]{\hat{\boldsymbol{v}}^\star\cdot \boldsymbol{n}^{\partial
          \hat{E}}}{\trace{\hat{w}}{\partial
          \hat{E}}{}}}{\norm[\sobh{1}{\hat{E}}]{\hat{w}}}
          \,.
      \end{split}
    \end{equation*}
    \fixBis{ Moreover, notice that the term
      $\sup_{\hat{w}\in\sobh{1}{\hat{E}}}\limits\frac{\scal[\partial
        \hat{E}]{\hat{\boldsymbol{v}}\cdot \boldsymbol{n}^{\partial
            \hat{E}}}{\trace{\hat{w}}{\partial
            \hat{E}}{}}}{\norm[\sobh{1}{\hat{E}}]{\hat{w}}}$ is a norm on
      $B^\partial(\Triang{E})$. Indeed, if
      $\sup_{\hat{w}\in\sobh{1}{\hat{E}}}\limits\frac{\scal[\partial
        \hat{E}]{\hat{\boldsymbol{v}}\cdot \boldsymbol{n}^{\partial
            \hat{E}}}{\trace{\hat{w}}{\partial
            \hat{E}}{}}}{\norm[\sobh{1}{\hat{E}}]{\hat{w}}}=0$ then
      $\hat{\boldsymbol{v}}\cdot \boldsymbol{n}^{\partial \hat{E}} = 0$ and
      $\hat{\boldsymbol{v}}=0$.}  Then, applying the above results to
    \cref{eq:vstarnnormestim_1}, \ourfixBis{recalling that
      $\sum_{\tau \in\Triang{E}}\norm[\lebl{\tau}]{\div \bs v}^2 = 0$ $\forall \bs v \in \mathrm{B}^{\partial}(\Triang{E})$}, using
    the equivalence of norms on finite dimensional spaces and a scaling
    argument, we get
    \begin{equation}
      \label{eq:vstarnormestim}
      \begin{split}
        \norm[\lebl{\partial E}]
        {\boldsymbol{v}^\star\cdot\boldsymbol{n}^{\partial E}}
          &\geq Ch_E^{-\frac12} \ourfixBis{h_E}
            \sup_{\hat{w}\in\sobh{1}{\hat{E}}}\frac{\scal[\partial
            \hat{E}]{\hat{\boldsymbol{v}}^\star\cdot \boldsymbol{n}^{\partial
            \hat{E}}}{\trace{\hat{w}}{\partial
            \hat{E}}{}}}{\norm[\sobh{1}{\hat{E}}]{\hat{w}}}
            \geq C h_E^{-\frac12}\ourfixBis{h_E}
            \norm[{\ourfixBis{\spaceV[\hat{E}]}}]{\hat{\boldsymbol{v}}^\star}
        \\
          &\ourfixBis{= C h_E^{-\frac12}  \left(\sum_{\tau\in\Triang{E}} \norm[\lebldouble{\tau}]{\boldsymbol{v}^\star}^2 + h^2_E\norm[\leblinfI]{\jmp[\intEdges]{\boldsymbol{v}^\star}}^2\right)^{\frac12}}
        \\
        & = C h_E^{-\frac12}\norm[{\ourfixBis{\spaceV}}]{\boldsymbol{v}^\star} \,,
      \end{split}
    \end{equation}
    where $C$ is independent of $h_E$ and of the choice of reference element by
    \cref{lem:compactnessSigma}.  The proof is thus concluded by applying the
    estimates \cref{eq:qnormestim} and \cref{eq:vstarnormestim} to
    \cref{eq:infsup-cont_firststep}.
  \end{proof}

  In the following, assuming \cref{eq:suffCondWellPos}, we prove the existence
  of an operator \ourfixBis{$\Pi_E$} satisfying \cref{eq:b-compatibility} and
  \cref{eq:Fortincontinuity}. First, we need some auxiliary results.
  \begin{definition}
    Let $\left\{ r_i \right\}_{i=1}^{\NVE-1}$ be a basis of $\RQE$.  Let us
    define the set of linear operators $D_i: \ourfixBis{\spaceV}\rightarrow \mathbb{R}$
    such that $\forall \boldsymbol{v}\in\ourfixBis{\spaceV}$
    \begin{equation*}
      D_i(\boldsymbol{v}):=\int_{\partial E} \limits \left(\boldsymbol{v}\cdot \boldsymbol{n}^{\partial E}\right)\trace{r_i}{\partial E}{} \,ds,
      \;\;\; \forall i=1,\ldots,\NVE-1 \,.
    \end{equation*}
  \end{definition}
\begin{lemma}
  If \ourfix{$(l+1)(l+2) - \dim\BadPoly{l}{E}\geq \NVE-1$}, there exists a set
  of functions $\boldsymbol{\pi}_j\in\PolyDouble{l}{E}$ defined by
  \begin{equation} \label{eq:indip functions for P}
    D_i(\boldsymbol{\pi}_j)=\delta_{ij}\;\; \forall i,j=1,\ldots,\NVE-1.
  \end{equation}
\end{lemma}
\begin{proof}
  Let $V^{M}_{l}(E)$ be the local mixed virtual element space of order $l$,
  defined in \cite{Beirao2016a}, i.e.
  \begin{equation*}
    \begin{gathered}
      V^{M}_{l}(E) := \{ \boldsymbol{v}\in \mathrm{H}(\mathrm{div};E)\cap \mathrm{H}(\mathrm{rot};E) : \trace{\boldsymbol{v}\cdot \boldsymbol{n}^e}{e}{} \in\Poly{l}{e} \forall e\in\E{E}, \\
      \mathrm{div}\boldsymbol{v}\in \Poly{l}{E} \mbox{ and }
      \mathrm{rot}\boldsymbol{v}\in\Poly{l-1}{E} \}.
    \end{gathered}
  \end{equation*}
  Notice that $\PolyDouble{l}{E}\subset V^{M}_{l}(E)$.  For each
  $\boldsymbol{v}\in V^{M}_{l}(E)$, the degrees of freedom of $\boldsymbol{v}$
  are defined \cite{Beirao2016a} by
  \begin{enumerate}
  \item
    $\int_e \boldsymbol{v}\cdot\boldsymbol{n}^e \,q \,\mbox{ds},\; \forall
    e\in\E{E},\, \forall q\in\Poly{l}{e}$,
  \item
    $\int_E \boldsymbol{v}\cdot \nabla p_l \,\mbox{dx},\; \forall p_l
    \in\Poly{l}{E}$,
  \item
    $\int_E\boldsymbol{v}\cdot \boldsymbol{p_l^{\perp}}\mbox{dx},\;
    \forall\boldsymbol{p_l^{\perp}}\in \{
    \boldsymbol{p_l^{\perp}}\in\PolyDouble{l}{E}: \int_E
    \boldsymbol{p_l^{\perp}} \cdot \nabla q\mbox{ dx}=0 \,\forall q
    \in\Poly{l+1}{E}\}$.
  \end{enumerate}
  The number of degrees of freedom defined by the first, the second and the
  third condition is, respectively, $(l+1)\NVE$, $\frac{(l+1)(l+2)}{2}-1$ and
  $\frac{(l-1)(l+2)}{2}+1$.  Globally,
  $\mathrm{dim}\, V^{M}_{l}(E) = (l+1)\NVE + l(l+2)$.

  Notice that a possible choice for the basis of
  $\mathcal{P}_{l}({\partial E}):= \{p\in\Poly{l}{e}, \forall e \in\E{E} \}$ is
  composed by the $\NVE -1 $ basis functions
  $\{\trace{r_i}{\partial E}{}\}_{i=1}^{\NVE-1}\subset\mathcal{Q}(\partial
  E)\subset\mathcal{P}_{l}({\partial E})$, completed by a choice of linearly
  independent functions
  $\{q^C_i\}_{i=\NVE}^{(l+1)\NVE}\subset\mathcal{P}_{l}({\partial E})$.  Hence,
  the first set of degrees of freedom can be split into two groups, i.e.
  \begin{itemize}
  \item
    $D_i(\boldsymbol{v}) = \int_{\partial E}
    \boldsymbol{v}\cdot\boldsymbol{n}^{\partial E} \,\trace{r_i}{\partial E}{}
    \,ds,\; \forall i=1,\ldots,\NVE -1$,
  \item
    $\int_{\partial E} \boldsymbol{v}\cdot\boldsymbol{n}^{\partial E} \,q^C_i
    \,\mbox{ds},\; \forall i=\NVE,\ldots,(l+1)\NVE$.
  \end{itemize}
  Let $j\in\{1,\ldots,\NVE-1\}$ and let $V^{R}(E;j) \subset V^{M}_{l}(E)$ be
  \begin{equation}
    V^{R}(E;j):= \{ \boldsymbol{v}\in V^{M}_{l}(E): D_i(\boldsymbol{v})=\ourfixBis{0}
    \; \forall i=1,\ldots, \NVE -1,\, \ourfixBis{i\neq j} \}.
  \end{equation}
  Notice that
  $\mathrm{dim}\, V^{R}(E;j) = \mathrm{dim}\, V^{M}_{l}(E) - (\NVE -1)
  \ourfixBis{+1} $.  Moreover, we define $V^{\perp\mathbb{P}_l}(E)$
  $\subset V^{M}_{l}(E)$, given by
  \begin{equation}
    V^{\perp\mathbb{P}_l}(E):= \{ \boldsymbol{v}\in V^{M}_{l}(E)\colon
    \mathrm{dof}(\boldsymbol{v})\cdot\mathrm{dof}(\boldsymbol{p}) =0 \,
    \forall \boldsymbol{p}\in\PolyDouble{l}{E}\ourfix{\setminus\BadPoly{l}{E}} \}
  \end{equation}
  where $\mathrm{dof}(\boldsymbol{v})$ denotes the vector of degrees of freedom
  of $\boldsymbol{v}\in V^{M}_{l}(E)$.  Notice that
  $\mathrm{dim}\,V^{\perp\mathbb{P}_l}(E) = \mathrm{dim}\, V^{M}_{l}(E)
  -\ourfix{\left((l+1)(l+2)-\dim\BadPoly{l}{E}\right)}$.  Since
  \ourfix{$(l+1)(l+2) - \dim\BadPoly{l}{E}\geq \NVE-1$}, then
  $\mathrm{dim}\, V^{R}(E;j) \ourfixBis{>}
  \mathrm{dim}\,V^{\perp\mathbb{P}_l}(E)$ and thus
  \begin{equation*} \label{eq:orth_system_w} \exists \boldsymbol{w}_j \in
    V^{R}(E;j) \cap
    \left(\PolyDouble{l}{E}\ourfix{\setminus\BadPoly{l}{E}}\right )
    \ourfixBis{,\, \boldsymbol{w}_j\neq \boldsymbol{0}}
    \,.
  \end{equation*}
  Then we can choose $\boldsymbol{\pi}_j=\boldsymbol{w}_j$ \ourfixBis{ such that $D_j(\boldsymbol{w}_j)=1$, this is possible since
    $D_j(\boldsymbol{w}_j)$ cannot be zero. Indeed, by contradiction let us
    suppose that $D_j(\boldsymbol{w}_j)=0$, then, by definition of
    $\BadPoly{\ell}{E}$ \cref{eq:defBadPolynomials},
    $\boldsymbol{w}_j\in\BadPoly{\ell}{E}$. This is a contradiction since
    $\boldsymbol{w}_j \in \PolyDouble{\ell}{E}\setminus \BadPoly{\ell}{E}$ and
    $\boldsymbol{w}_j\neq \boldsymbol{0}$. }
\end{proof}

In the following proposition we provide a definition of \ourfixBis{$\Pi_E$ and
prove \cref{eq:b-compatibility} and \cref{eq:Fortincontinuity}}.
\begin{proposition}\label{prop:def Pi2}
  Under the hypothesis of \cref{PiZeroEnhInjectivity}, let us define
  $\ourfixBis{\Pi_E}:\ourfixBis{\spaceV}\rightarrow\PolyDouble{l}{E}$ such that
  $\forall \boldsymbol{v}\in\ourfixBis{\spaceV}$
  \begin{equation*}
    \ourfixBis{\Pi_E}\boldsymbol{v} := \sum_{i=1}^{\NVE-1}  D_i(\boldsymbol{v})\boldsymbol{\pi}_i \, ,
  \end{equation*}
  where $\boldsymbol\pi_i$ satisfy \cref{eq:indip functions for P}.  Then
  $\ourfixBis{\Pi_E}$ satisfies \ourfixBis{\cref{eq:b-compatibility} and
    \cref{eq:Fortincontinuity}}.
\end{proposition}

\begin{proof}
  Since
  \begin{equation}
    \label{eq:Di_DPi}
    \forall \boldsymbol{v}\in \ourfixBis{\spaceV}, \;\;
    D_i(\ourfixBis{\Pi_E}\boldsymbol v)=D_i(\boldsymbol{v}) \;\; \forall i=1,\ldots,\NVE-1,
  \end{equation}
  let us check that $\ourfixBis{\Pi_E}$ satisfies \cref{eq:b-compatibility},
  indeed by construction
  $\forall r_i\in\RQE, \,i=1,\ldots,\NVE-1, \,\forall\boldsymbol{v}\in
  \ourfixBis{\spaceV}$:
  \begin{align*}
    b(r_i,\ourfixBis{\Pi_E}\boldsymbol{v}-\boldsymbol{v})
    &=
      \int_{\partial E} \limits r_i\,\left(\ourfixBis{\Pi_E}\boldsymbol{v}-\boldsymbol{v}\right)\cdot n^{\partial E} \, dx
      =D_i(\ourfixBis{\Pi_E}\boldsymbol{v}-\boldsymbol{v})
      =0.
  \end{align*}
  Furthermore, let us consider
  $\widehat{\ourfixBis{\Pi_E}\boldsymbol{v}} \ourfixBis{ = \Pi_E\boldsymbol{v}
    \circ F}$ defined on the reference polygon $\hat{E}$.  Applying the
  linearity of the definition of the mapping $F:\hat{E}\rightarrow E$, presented
  in \cref{def:affinemapping}, we have
  \begin{equation} \label{eq:PiE_onRef}
    \widehat{\ourfixBis{\Pi_E}\boldsymbol{v}} =\left(\sum_{i=1}^{\NVE-1}\limits
      D_i\left(\boldsymbol{v}\right)\boldsymbol{\pi}_i\right)\circ F
    =\sum_{i=1}^{\NVE-1}\limits
    D_i\left(\boldsymbol{v}\right)\left(\boldsymbol{\pi}_i\circ F\right)
    \ourfixBis{= \sum_{i=1}^{\NVE-1}\limits
      D_i\left(\boldsymbol{v}\right)\hat{\boldsymbol{\pi}}_i} \,.
  \end{equation}
  \ourfixBis{ Then, applying \cref{lem: continuity of b}, we have
    $\forall i = 1, \ldots, \NVE -1$
    \begin{equation} \label{eq:continuity of Di} \abs{D_i(\bs v)} =
      \abs{\int_{\partial E} \limits \left(\boldsymbol{v}\cdot
          \boldsymbol{n}^{\partial E}\right)\trace{r_i}{\partial E}{} \,ds } = h_E
      \abs{b(\hat{r}_i,\hat{\bs{v}})} \leq C_b h_E\norm[\HT{\hat{E}}]{\hat{r}_i}
      \norm[{\spaceV[\hat{E}]}]{\hat{\boldsymbol{v}}}.
    \end{equation}
    Then, we want to prove the continuity of $\Pi_E\boldsymbol{v}$.  Since
    \begin{equation*}
      \Pi_E\boldsymbol{v}\in [\Poly{\ell}{E}]^2\imply \Pi_E\boldsymbol{v}\in\cont{E}\imply\norm[\leblinfI]{\jmp[\intEdges]{\Pi_E\boldsymbol v}}=0,
    \end{equation*}
    applying \cref{eq:PiE_onRef} and \cref{eq:continuity of Di}, we have
    \begin{equation}\label{eq:contPiE_1}
      \begin{aligned}
        \norm[\spaceV]{\Pi_E\boldsymbol{v}} ^2
        &= \norm[\lebldouble{E}]{\Pi_E\boldsymbol{v}}^2+
          \norm[\lebl{E}]{\div \Pi_E\boldsymbol{v}}^2
        \\
        &= h^2_E \norm[\lebl{\hat{E}}]{\sum_{i=1}^{\NVE-1}\limits D_i\left(\boldsymbol{v}\right)\hat{\boldsymbol{\pi}}_i}^2
          +
          \norm[\lebl{\hat{E}}]{\hat{\nabla}\cdot\left(\sum_{i=1}^{\NVE-1}\limits D_i\left(\boldsymbol{v}\right)\hat{\boldsymbol{\pi}}_i\right)}^2
        \\
        &\leq C\sum_{i=1}^{\NVE-1}\limits \abs{D_i\left(\boldsymbol{v}\right)}^2\left(
          h^2_E\norm[\lebldouble{\hat{E}}]{\hat{\bs{\pi}}_i}^2 + 
          \norm[\lebl{\hat{E}}]{\hat{\nabla}\cdot\hat{\bs{\pi}}_i}^2
          \right)
        \\
        &\leq C \NVmax \max_{i=1,\ldots,\NVE-1}\limits \left\{ \norm[\HT{\hat{E}}]{\hat{r}_i}^2
          \norm[{\spaceV[\hat{E}]}]{\hat{\bs{\pi}}_i}^2
          \right\}   h_E^2 \norm[{\spaceV[\hat{E}]}]{\hat{\boldsymbol{v}}} \,.
      \end{aligned}
    \end{equation}
  } We set
  $ C(\hat{E}) := \max_i \limits \norm[\HT{\hat{E}}]{\hat{r}_i}
  \max_i\limits\norm[\ourfixBis{\spaceV[\hat{E}]}]{\widehat{\boldsymbol{\pi}}_i}$. This
  is a continuous function on the set of admissible reference elements $\Sigma$,
  which is a compact set by \cref{lem:compactnessSigma}.  Indeed,
  $\norm[\HT{\hat{E}}]{\hat{r}_i}$ is a continuous function
  $\forall i =1,\ldots, \NVE -1$ on $\Sigma$.  Moreover, by definition,
  $\widehat{\boldsymbol{\pi}}_i$ depends continuously on the set
  $\{\hat{r}_i\}_{i=1}^{\NVE-1}$.  Then there exists
  $M = \max_{\hat{E}\in\Sigma} C(\hat{E}) >0$.  \ourfixBis{Finally, starting
    from \cref{eq:contPiE_1}, it results $\exists C>0$ such that
    \begin{equation}
      \begin{aligned}
        \norm[\spaceV]{\Pi_E\boldsymbol{v}} ^2
        &\leq 
          C h_E^2 \norm[{\spaceV[\hat{E}]}]{\hat{\boldsymbol{v}}}
        \\
        &\leq C  \left(h_E^2\sum_{\hat{\tau}\in\Triang{\hat{E}}} \norm[\lebldouble{\hat{\tau}}]{\hat{\bs{v}}}^2 +
          \sum_{\hat{\tau}\in\Triang{\hat{E}}} \norm[\lebl{\hat{\tau}}]{\hat{\nabla}\cdot\hat{\bs{v}}}^2+ h^2_E\norm[\leblinfIReferenceElement]{\jmp[\intEdgesReferenceElement]{\hat{\bs{v}}}}^2\right)\\
        & = C\norm[\spaceV]{\bs v}^2\,.
      \end{aligned}
    \end{equation}
  }
\end{proof}




\ourfix{%
  \subsection{Upper bound on the dimension of $\BadPoly{l}{E}$}
  \label{sec:upperboundbadpoly}
  \begin{theorem}
    \label{th:badpolyupperbound}
    Let $E\in\Mh$ and $l\in\mathbb{N}$. Then,
    \begin{equation}
      \label{eq:upperbounddimBadPoly}
      \dim \BadPoly{l}{E} \leq l(l+1) \,.
    \end{equation}
  \end{theorem}
  \begin{proof}
    Consider the following space of harmonic polynomials:
    \begin{equation*}
      \HPoly{l+1}{E} = \{ p \in \Poly{l+1}{E}\colon \Delta p = 0 \} \,.
    \end{equation*}
    It is known that $\dim \HPoly{l+1}{E} = 2l+3$, thus the space of its
    gradients $\nabla \HPoly{l+1}{E}$ satisfies
    $\dim\nabla \HPoly{l+1}{E} = 2l+2$. We prove the thesis by showing that
    $\nabla \mathbb{H}_{l+1}\left(E\right) \cap \BadPoly{l}{E} =
    \{\mathbf{0}\}$. The result will follow by difference, since
    \begin{equation*}
      \dim \PolyDouble{l}{E} - \dim\nabla \mathbb{H}_{l+1}\left(E\right) =
      l(l+1) \,.
    \end{equation*}
    Then, let $q_{l+1} \in\HPoly{l+1}{E}$ and suppose
    $\nabla q_{l+1} \in \BadPoly{l}{E}$. Then, $\forall \varphi \in \sobh{1}{E}$
    such that $\trace{\varphi}{\partial E}{} \in \mathcal{Q} (\partial E)$ it
    holds
    \begin{equation}
      \label{eq:gradHarmonicGradPhi}
      0  = \scal[E]{-\Delta q_{l+1}}{\varphi} =
      \scal[E]{\nabla q_{l+1}}{\nabla \varphi} -
      \int_{\partial E}\nabla q_{l+1}\cdot \boldsymbol{n}^{\partial E}
      \trace{\varphi}{\partial E}{} = \scal[E]{\nabla q_{l+1}}{\nabla \varphi}\,.
    \end{equation}
    Then, let $\omega \subset E$ be any open subset of $E$ such that
    $\partial \omega \cap \partial E = \emptyset$ and let
    \fixBis{$c = \frac{1}{\abs{E}}\int_E q_{l+1}$}. Let $b_\omega\in\sobh{1}{E}$ be
    the bubble function defined as:
    \begin{equation*}
      \begin{cases}
        -\Delta b_\omega = q_{l+1}-c & \text{in $\omega$}\,,
        \\
        b_\omega = 0 & \text{on $E\setminus \omega$} \,.
      \end{cases}
    \end{equation*}
    For any $\chi \in \mathcal{Q} (\partial E)$, let
    $\varphi_1,\,\varphi_2\in\sobh{1}{E}$ be two functions such that
    $\trace{\varphi_1}{\partial E}{}=\trace{\varphi_2}{\partial E}{} = \chi$ and
    $\eval{(\varphi_1 - \varphi_2)}{E} = b_\omega$. Then, by
    \cref{eq:gradHarmonicGradPhi},
    \begin{equation*}
      \scal[E]{\nabla q_{l+1}}{\nabla b_\omega} =
      \scal[E]{\nabla q_{l+1}}{\nabla (\varphi_1-\varphi_2)} = 0 \,.
    \end{equation*}
    By the definition of $b_\omega$, we get
    \begin{equation*}
      \norm[\lebl{\omega}]{q_{l+1}-c}^2 = \scal[\omega]{q_{l+1}-c}{q_{l+1}-c}
      =  -\scal[\omega]{q_{l+1} - c}{\Delta b_\omega}
      =  \scal[E]{\nabla q_{l+1}}{\nabla b_\omega} = 0\,.
    \end{equation*}
    which means that $\eval{q_{l+1}}{\omega} = c$. In particular, $q_{l+1}=c$ in
    at least $\fixBis{\frac{1}{2}}(l+2)(l+3) = \dim\Poly{l+1}{E}$ distinct points in the interior of
    $E$. Thus, $q_{l+1}=c$ and $\nabla q_{l+1} = \mathbf{0}$.
  \end{proof}
}%
\ourfix{
  \begin{remark}
    If we consider $E^\star$ to be the hexagon having vertices
    \begin{equation*}
      x_i = \begin{pmatrix} \cos\left(\frac{(i-1)\pi}{6}\right)
        &\sin\left(\frac{(i-1)\pi}{6}\right)\end{pmatrix} \quad,
      i\in\{1,\ldots,6\}\,,
    \end{equation*}
    an explicit computation yields
    \begin{equation*}
      \BadPoly{1}{E^\star} = \mathrm{span}\left\{
        \begin{pmatrix} x\\y \end{pmatrix} , \begin{pmatrix} y\\-x \end{pmatrix}
      \right\} \,,
    \end{equation*}
    which implies that the estimate provided by \cref{th:badpolyupperbound} is
    optimal for $l=1$.
  \end{remark}
}
\ourfixBis{
  \begin{remark}
  Notice that the upper bound on the dimension of $\BadPoly{l}{E}$ given by $l(l+1)$ in \cref{eq:upperbounddimBadPoly} is independent of the geometry of the polygon. Then,
  it allows us to obtain the sufficient condition $2l+2\geq\NVE - 1 $ for \cref{PiZeroEnhInjectivity}, which is robust with respect to the element shape.
  \end{remark}
  }

\subsection{Coercivity of the discrete bilinear form}
\label{sec:coercivity}

In this section we prove the coercivity of the discrete problem defined by \cref{eq:discrVarForm}
with respect to the standard $\sobh[0]{1}{\Omega}$ norm, denoted by
\begin{equation*}
  \norm[{\sobh[0]{1}{\Omega}}]{V} = \norm[\lebldouble{\Omega}]{\nabla V} \quad \forall V \in \sobh[0]{1}{\Omega}\,.
\end{equation*}
Let
\begin{equation*}
  \norm[\boldell]{v} := \left(
    \sum_{E\in\Mh} \norm[\lebldouble{E}]{\proj{\boldell(E)}{E} \nabla v}^2
  \right)^{\frac12}
  \quad \forall v \in \V{1,\boldell} \,.
\end{equation*}
We have the following result.
\begin{proposition}
  Suppose $\boldell$ satisfies \cref{eq:suffCondWellPos} $\forall E \in \Mh$.  Then,
  $\norm[\boldell]{\cdot}$ is a norm on $\V{1,\boldell}$.
\end{proposition}
\begin{proof}
  Let $v\in\V{1,\boldell}$ be given. It is clear from its definition that $\norm[\boldell]{v}$ is a
  semi-norm.  Applying \cref{PiZeroEnhInjectivity} and since $v\in \sobh[0]{1}{\Omega}$, we
  have that
  \begin{equation*}
    \norm[\boldell]{v} = 0 \imply \norm[{\sobh[0]{1}{\Omega}}]{v} = 0 \imply v = 0\,.
  \end{equation*}
  ~
\end{proof}

\begin{lemma}
  \label{lem:ahBounds}
  We have that
  \begin{equation}
  \label{eq:upperBoundBoldellNorm}
  \norm[\boldell]{v} \leq \norm[{\sobh[0]{1}{\Omega}}]{v} \quad \forall v \in \V{1,\boldell} \,.
\end{equation}
Moreover, if $\boldell(E)$ satisfies \cref{eq:suffCondWellPos}
$\forall E \in \Mh$, then
  \begin{equation}
    \label{eq:lowerBoundBoldellNorm}
  \exists c_\ast > 0 \colon
  \norm[\boldell]{v} \geq c_\ast \norm[{\sobh[0]{1}{\Omega}}]{v}
  \quad \forall v \in \V{1,\boldell}\,,
\end{equation}
where $c_\ast$ does not depend on $h$.
\end{lemma}
\begin{proof}
  Relation \cref{eq:upperBoundBoldellNorm} follows immediately by the definition of $\proj{l}{E}$
  and an application of the Cauchy-Schwarz inequality. Indeed, let $E\in\Mh$, then
  \begin{equation*}
    \norm[E]{\proj{l}{E} \nabla v}^2 \! = \!\scal[E]{\proj{l}{E} \nabla v}{\proj{l}{E} \nabla v} \!=\!
    \scal[E]{\nabla v}{\proj{l}{E}\nabla v}
    \leq \norm[\lebldouble{E}]{\nabla v} \norm[\lebldouble{E}]{\proj{l}{E}\nabla v} \,.
  \end{equation*}
  On the other hand, by standard scaling arguments we have
  \begin{equation*}
    \norm[\boldell]{v}^ 2 =
    \sum_{E\in\Mh} \norm[\lebldouble{E}]{\proj{l}{E}\nabla v}^2
    = \sum_{E\in\Mh} \norm[\lebldouble{\hat{E}}]{\projhat{l}{\hat{E}}\hat{\nabla}\left(\hat{v}- P_0 (\hat{v})\right)}^2.
  \end{equation*}
  Notice that $\forall \hat{E}\in\Sigma$, where $\Sigma$ is the set of admissible reference
  elements, $\hat{v}- \mathrm{P}_0 (\hat{v})\in\V[\hat{E},\mathrm{P}_0]{1,l}$.  Moreover,
  $\forall \hat{w}\in\V[\hat{E},\mathrm{P}_0]{1,l}$ both
  $\norm[\lebldouble{\hat{E}}]{\projhat{l}{\hat{E}}\hat{\nabla}\hat{w}}$ and
  $\norm[\lebldouble{\hat{E}}]{\hat{\nabla} \hat{w}}$ are norms.  Then we can apply, by standard
  arguments about the equivalence of norms on finite dimensional spaces, we obtain
  $\forall \hat{E}\in\Sigma$
  \begin{equation}
    \norm[\lebldouble{\hat{E}}]{\projhat{l}{\hat{E}}\hat{\nabla}\hat{w}}
    \geq C(\hat{E})\norm[\lebldouble{\hat{E}}]{\hat{\nabla} \hat{w}}
  \end{equation}
  where
  \begin{equation}
    C(\hat{E})=
    \frac{\min_{\hat{z}\in\V[\hat{E},\mathrm{P}_0]{1,l} : \norm[\ourfix{l^2}]{\dof{}{\hat{z}}}=1}\norm[\lebldouble{\hat{E}}]{\projhat{l}{\hat{E}}\hat{\nabla}\hat{z}}}
    {\sqrt{\NVE -1}\max_{i=1,\ldots,\NVE -1} \norm[\lebldouble{\hat{E}}]{\hat{\nabla} \hat{\psi}_i}}\, .
  \end{equation}
  $C(\hat{E})$ is a continuous function on $\Sigma$, which is a compact set by \cref{lem:compactnessSigma}.  Indeed, $\projhat{l}{\hat{E}}$ is continuous on $\Sigma$, as well as
  functions in $\V[\hat{E},\mathrm{P}_0]{1,l}$ following proofs of \cite[Lemma 4.9]{Cangiani2016}
  and \cite[Lemma 4.5]{Beirao2015}.  Moreover, $C(\hat{E})>0$, $\forall \hat{E}\in\Sigma$.  Indeed,
  applying \cref{prop: Brezzi _1}, it holds that
  $\forall \hat{z}\in \V[\hat{E},\mathrm{P}_0]{1,l} : \norm[\ourfix{l^2}]{\dof{}{\hat{z}}}=1$,
  \begin{equation*}
    \begin{aligned}
    \norm[\lebldouble{\hat{E}}]{\projhat{l}{\hat{E}}\hat{\nabla}\hat{z}}^2
   &= \scal[\hat{E}]{\hat{\nabla}\hat{z}}{\projhat{l}{\hat{E}}\hat{\nabla}\hat{z}}
   = \scal[\partial \hat{E}]{\hat{z}}{\projhat{l}{\hat{E}}\hat{\nabla}\hat{z} \cdot \boldsymbol{n}^{\partial\hat{E}}}
=b(\hat{z}_R,\projhat{l}{\hat{E}}\hat{\nabla}\hat{z})
   \\
   &\geq \beta \norm[\lebldouble{\hat{E}}]{\projhat{l}{\hat{E}}\hat{\nabla}\hat{z}} \norm[\HT{\hat{E}}]{\hat{z}_R}>0 \,,
 \end{aligned}
\end{equation*}
  where $\hat{z}_R$ is the lifting of $\trace{\hat{z}}{\partial \hat{E}}{}$ on $\RQE[\hat{E}]$.
  Then, $\exists m>0$ such that $m := \min_{\hat{E}\in\Sigma} C(\hat{E})$.
  Finally, by standard scaling argument we obtain
  \begin{equation}
    \norm[\boldell]{v}^ 2 \geq m^2 \sum_{E\in\Mh} \norm[\lebldouble{\hat{E}}]{\hat{\nabla}\left(\hat{v}-\mathrm{P}_0 (\hat{v})\right)}^2
    = m^2 \norm[{\sobh[0]{1}{\Omega}}]{v} \,.
  \end{equation}
  ~
\end{proof}

In the following theorem, we provide a proof of the continuity and the coercivity of the discrete bilinear form. The coercivity property follows from \cref{lem:ahBounds}.

\begin{theorem}
  Let $\a[h]{}{}$ be the bilinear form defined by \cref{eq:defah}. Then,
  \begin{equation}
    \label{eq:ahContinuity}
      \a[h]{w}{v} \leq  \norm[{\sobh[0]{1}{\Omega}}]{w}\norm[{\sobh[0]{1}{\Omega}}]{v}
    \forall w,v \in\V{1,\boldell}\, .
  \end{equation}
  Moreover, suppose $\boldell(E)$ satisfies \cref{eq:suffCondWellPos}
  $\forall E \in \Mh$. Then,
  \begin{equation}
    \label{eq:ahCoercivity}
    \exists C>0, \mbox{ independent of } h \colon \a[h]{w}{w} \geq C \norm[{\sobh[0]{1}{\Omega}}]{w}^2
    \forall w \in\V{1,\boldell}\,.
  \end{equation}
\end{theorem}
\begin{proof}
  Let $w,v \in\V{1,\boldell}$ be given. Applying the Cauchy-Schwarz inequality and
  \cref{eq:upperBoundBoldellNorm} we get
  \begin{equation*}
    \begin{split}
      \a[h]{w}{v} &= \sum_{E\in\Mh} \scal[E]{\proj{\boldell(E)}{E} \nabla w}{\proj{\boldell(E)}{E}
        \nabla v}
      \\
      &\leq \sum_{E\in\Mh} \norm[\lebldouble{E}]{\proj{\boldell(E)}{E} \nabla w} \norm[\lebldouble{E}]{\proj{\boldell(E)}{E}
        \nabla v}
      \\
      &\leq \norm[\boldell]{w} \norm[\boldell]{v} \leq \norm[{\sobh[0]{1}{\Omega}}]{w}
      \norm[{\sobh[0]{1}{\Omega}}]{v} \,.
    \end{split}
  \end{equation*}
  Moreover, assuming that $\boldell(E)$ satisfies \cref{eq:suffCondWellPos}
  $\forall E \in \Mh$, we can apply the lower bound in
  \cref{eq:lowerBoundBoldellNorm} and get
  \begin{equation*}
    \a[h]{w}{w} = \norm[\boldell]{w}^2 \geq
    \left( c_\ast \right)^2 \norm[{\sobh[0]{1}{\Omega}}]{w}^2\,.
  \end{equation*}
  ~
\end{proof}
This theorem implies that the bilinear form $a_h$ of the problem
\cref{eq:discrVarForm} satisfies the hypothesis of \fix{the} Lax-Milgram
theorem, \fix{hence} the problem admits a unique solution.

\section{A priori error estimates}
\label{sec:apriori}
In this section we derive error estimates for the proposed method, in $\sobh[0]{1}{}$ norm and in
the standard $\mathrm{L}^2$ norm.  First, we recall classical results for Virtual Element Methods
concerning the interpolation error and the polynomial projection error (see
\cite{Cangiani2017,Beirao2015b}). 
\begin{lemma}\label{lem:interpolationError}
  Let \ourfix{$U\in\sobh{2}{\Omega}$}, then there exists $C>0$
  such that $\forall \,h$, $\exists\UI\in\V{1,\boldell}$ \fix{satisfying}
  \begin{equation}\label{eq:interpolationError}
    \norm[\lebl{\Omega}]{U-\UI}+h\norm[{\sobh[0]{1}{\Omega}}]{U-\UI} \leq C h^2\seminorm[2]{U}.
  \end{equation}
\end{lemma}
\begin{proof}
  The proof of this result can be obtained following the same arguments as in \cite[Theorem
  11]{Cangiani2017}.
\end{proof}
\begin{lemma}[{\cite[Lemma 5.1]{Beirao2015b}}]
  Let \ourfix{$U\in\sobh{2}{\Omega}$}, then there exist
  $C_1,C_2>0$ such that
  \begin{align}\label{eq:polynomialErrorGrad}
    \norm[\lebl{\Omega}]{\proj{\boldell}{} \nabla U-\nabla U}&\leq C_1 h\seminorm[2]{U},
    \\
    \label{eq:polynomialErrorPiZero}
    \norm[\lebl{\Omega}]{\proj{0}{} U- U}&\leq C_2 h\norm[{\sobh[0]{1}{\Omega}}]{U}.
  \end{align}
\end{lemma}


\begin{theorem} \label{thm:H1estimate}
Let $U\in\sobh{2}{\Omega}\cap\sobh[0]{1}{\Omega}$ and $f\in\lebl{\Omega}$ be the solution and the right-hand side of \cref{eq:modelProblem}, respectively. \fix{Then,} $\exists C>0$ such that the unique solution $u\in\V{1,\boldell}$ to problem \cref{eq:discrVarForm} satisfies the following error estimate:
\begin{equation}\label{eq:H1estimate}
\norm[{\sobh[0]{1}{\Omega}}]{U-u}\leq Ch\left(\seminorm[2]{U}+\norm[\lebl{\Omega}]{f}\right) \,.
\end{equation}
\end{theorem}
\begin{proof}
Let $\UI$ be given by \cref{lem:interpolationError}. Applying the triangle inequality, we have
\begin{equation}\label{eq:firstStep-triangleInequality}
\norm[{\sobh[0]{1}{\Omega}}]{U-u}\leq \norm[{\sobh[0]{1}{\Omega}}]{U-\UI} + \norm[{\sobh[0]{1}{\Omega}}]{\UI-u}.
\end{equation}
We deal with the two terms separately. The first one can be bounded applying \cref{eq:interpolationError}, i.e.
\begin{equation}\label{eq:firstTermofProof}
\norm[{\sobh[0]{1}{\Omega}}]{U-\UI} \leq Ch\seminorm[2]{U}.
\end{equation}
On the other hand, in order to deal with the second term of \cref{eq:firstStep-triangleInequality} let $\varepsilon=\UI-u$. First, applying the coercivity of the bilinear form $a_h$ \cref{eq:ahCoercivity} and the discrete problem \cref{eq:discrVarForm}, we have that $\exists C>0$:
\begin{equation} \label{eq:zeroStep-EpsilonEstimate}
C\norm[{\sobh[0]{1}{\Omega}}]{\varepsilon}^2\leq \ah{\varepsilon}{\varepsilon}=\ah{\UI}{\varepsilon}-\ah{u}{\varepsilon} = \ah{\UI}{\varepsilon} -\sum_{E\in\Mh}\scal[E]{f}{\proj{0}{E}\varepsilon}.
\end{equation}
Applying the definition of the $\mathrm{L}^2$ projectors and adding and subtracting terms, i.e. $\proj{l}{E} \nabla U$ and $\nabla U$, we have
\begin{align*}
\ah{\varepsilon}{\varepsilon} &= \ah{\UI-U}{\varepsilon} + \ah{U}{\varepsilon} -\sum_{E\in\Mh}\scal[E]{\proj{0}{E}f}{\varepsilon} \\
&=\ah{\UI-U}{\varepsilon} + \sum_{E\in\Mh} \scal[E]{\proj{l}{E}\nabla U-\nabla U}{\nabla \varepsilon} 
+\scal[E]{\nabla U}{\nabla \varepsilon} - \scal[E]{\proj{0}{E}f}{\varepsilon} \\
&=\ah{\UI-U}{\varepsilon} + \sum_{E\in\Mh} \scal[E]{\proj{l}{E}\nabla U-\nabla U}{\nabla \varepsilon} +\scal[E]{f-\proj{0}{E}f}{\varepsilon}.
\end{align*}
Let us consider the last three terms separately.
The first one can be bounded applying \cref{eq:ahContinuity} and \cref{eq:interpolationError}, i.e.
\begin{equation} \label{eq:firstTerm-EpsilonEstimate}
\ah{\UI-U}{\varepsilon} \leq C \norm[{\sobh[0]{1}{\Omega}}]{\UI-U}\norm[{\sobh[0]{1}{\Omega}}]{\varepsilon} \leq Ch\seminorm[2]{U}\norm[{\sobh[0]{1}{\Omega}}]{\varepsilon}.
\end{equation}
Applying the Cauchy-Schwarz inequality and \cref{eq:polynomialErrorGrad}, the second term can be bounded as follows:
\begin{equation}\label{eq:secondTerm-EpsilonEstimate}
\begin{aligned}
\sum_{E\in\Mh} \scal[E]{\proj{l}{E}\nabla U-\nabla U}{\nabla \varepsilon} 
&\leq \sum_{E\in\Mh} \norm[\lebl{E}]{\proj{l}{E}\nabla U-\nabla U} \norm[{\sobh[0]{1}{E}}]{\varepsilon} \\
&\leq Ch\seminorm[2]{U}\norm[{\sobh[0]{1}{\Omega}}]{\varepsilon} \,.
\end{aligned}
\end{equation}
The last term can be bounded applying the definition of $\proj{0}{E}$, the Cauchy-Schwarz inequality and \cref{eq:polynomialErrorPiZero}, i.e.
\begin{align}\notag
\sum_{E\in\Mh}\scal[E]{f-\proj{0}{E}f}{\varepsilon} &= \sum_{E\in\Mh}\scal[E]{f}{\varepsilon-\proj{0}{E}\varepsilon} \\ \label{eq:thirdTerm-EpsilonEstimate}
&\leq \sum_{E\in\Mh} \norm[\lebl{E}]{f} \norm[\lebl{E}]{\varepsilon-\proj{0}{E}\varepsilon} \leq Ch\norm[\lebl{\Omega}]{f}\norm[{\sobh[0]{1}{\Omega}}]{\varepsilon} \,.
\end{align}
Finally, applying together \cref{eq:firstTerm-EpsilonEstimate},\cref{eq:secondTerm-EpsilonEstimate} and \cref{eq:thirdTerm-EpsilonEstimate} into \cref{eq:zeroStep-EpsilonEstimate} and simplifying, we have
\begin{equation}\label{eq:secondTermofProof}
\norm[{\sobh[0]{1}{\Omega}}]{\varepsilon} \leq C h \left(\seminorm[2]{U}+ \norm[\lebl{\Omega}]{f}\right).
\end{equation}
Considering together \cref{eq:firstTermofProof} and \cref{eq:secondTermofProof} we prove \cref{eq:H1estimate}.
\end{proof}
\begin{theorem} \label{thm:L2estimate} Let \ourfix{$U\in\sobh{2}{\Omega}\cap\sobh[0]{1}{\Omega}$}
  and $f\in\sobh{1}{\Omega}$ be the solution and the right-hand side of \cref{eq:modelProblem},
  respectively. \fix{Then}, $\exists C>0$ such that the unique solution $u\in\V{1,\boldell}$ to
  problem \cref{eq:discrVarForm} satisfies the following error estimate:
\begin{equation}\label{eq:L2estimate}
\norm[\lebl{\Omega}]{U-u}\leq Ch^2\left(\seminorm[2]{U}+\norm[{\sobh[0]{1}{\Omega}}]{f}\right) \,.
\end{equation}
\end{theorem}
\begin{proof}
The proof of this theorem can be found in the supplementary materials of this paper.
\end{proof}
\begin{remark}
Denoting by $\proj{1}{E}$ the $\lebl{}$-projector from $\lebl{E}$ to $\Poly{1}{E}$, we can define the discrete problem \cref{eq:discrVarForm} as
\begin{equation*}
  \ah{u}{v} = \sum_{E\in\Mh}\scal[E]{f}{\proj{1}{E}v}
  \quad \forall v \in \V{1,\boldell}\,,
\end{equation*}
and we can require $f\in\lebl{\Omega}$ so \cref{eq:L2estimate} still holds as 
\begin{equation*}
\norm[\lebl{\Omega}]{U-u}\leq Ch^2\left(\seminorm[2]{U}+\norm[{\lebl{\Omega}}]{f}\right) \,.
\end{equation*}
\end{remark}

\begin{remark}[Extension to more general elliptic problems]
  \label{sec:diffreact}
  Consider the following diffusion-reaction model:
  \begin{equation}
    \label{eq:diffreact}
    \begin{cases}
      -\Delta U + U = f & \text{in $\Omega$}\,,
      \\
      U = 0 & \text{on $\partial \Omega$}\,.
    \end{cases}
  \end{equation}
  The coercivity of the bilinear form defined by \cref{eq:defahE} and \cref{eq:defah} allows us to discretize it as: find $u \in \V{1,\boldell}$ such that
  \begin{equation}
    \label{eq:discr-diffreact}
    \a[h]{u}{v} + \sum_{E\in\Mh}\scal[E]{\proj{0}{E}u}{\proj{0}{E}v} =
    \scal[E]{f}{\proj{0}{E}v} \quad \forall v \in \V{1,\boldell} \,.
  \end{equation}
  If $\boldell$ satisfies \cref{eq:suffCondWellPos} locally on each polygon, we can prove the
  well-posedness of \cref{eq:discr-diffreact} following \cite[Lemma
  5.7]{Beirao2015b}. Optimal order a priori error estimates can be proved as in
  \cite[Theorem 5.1 and 5.2]{Beirao2015b}, using the interpolation result given by \cref{lem:interpolationError}. In \cref{sec:numres:diffreact} we
  assess numerically the validity of such results.
\end{remark}

\section{Numerical Results}
\label{sec:numericalresults}

\ourfix{This section is devoted to assess the theoretical results reported
  previously. First, we consider single polygons and investigate numerically
  which is the minimum degree $l$ providing coercivity, then we carry out some
  convergence tests.}
\subsection{Coercivity tests}
\label{sec:numres:coercivity}
\begin{table}
  \centering
  \caption{\ourfix{Sufficient $l$ for regular polygons up to 20 edges}}
  \ourfix{ \begin{tabular}{c|cccccccccc}
  $\NVE$ & $3$ &$4, 5$ &$6, 7$ &$8, 9$ &$10, 11$ &$12, 13$ &$14, 15$ &$16, 17$ &$18, 19$ &$20$
  \\
  $\hat{\ell}(\NVE)$ & $0$ & $1$ & $2$  & $3$ & $4$ & $5$ & $6$ & $7$ & $8$ & $9$
  \\
  $l$ & $0$ &$1$ &$2$  &$3$ &$4$ &$5$ &$6$  &$7$ &$8$ &$9$
  \\
  $\check{\ell}(\NVE)$ & $0$ & $1$ & $1$ & $2$ & $2$ & $2$ & $3$ & $3$ & $3$ & $3$
\end{tabular}
 }
  \label{tab:ell-regularPoly}%
\end{table}

\begin{table}
  \centering
  \caption{\ourfix{Sufficient $l$ for non-regular convex polygons up to 20
      edges}}
  \ourfix{ \begin{tabular}{c|cccccccccc}
  $\NVE$ & $3$ &$4, 5$ &$6, 7$ &$8, 9$ &$10, 11$ &$12, 13$ &$14, 15$ &$16, 17$ &$18, 19$ &$20$
  \\
  $\hat{\ell}(\NVE)$ & $0$ & $1$ & $2$  & $3$ & $4$ & $5$ & $6$ & $7$ & $8$ & $9$
  \\
  $l$ & $0$ & $1$ & $1$ & $2$ & $2$ & $2$ & $3$ & $3$ & $3$ & $3$
  \\
  $\check{\ell}(\NVE)$ & $0$ & $1$ & $1$ & $2$ & $2$ & $2$ & $3$ & $3$ & $3$ & $3$
\end{tabular}
 }
  \label{tab:ell-nonregularPoly}
\end{table}

\begin{table}
  \centering
  \caption{\ourfix{Sufficient $l$ for polygons with aligned edges up to 12
      edges}}
  \ourfix{ \begin{tabular}{c|cccccccccc}
  $\NVE$ & $3$ &$4, 5$ &$6, 7$ &$8, 9$ &$10, 11$ & $12$ 
  \\
  $\hat{\ell}(\NVE)$ & $0$ & $1$ & $2$  & $3$ & $4$ & $5$
  \\
  $l$ & $0$ & $1$ & $2$ & $2$ & $3$ & $4$
  \\
  $\check{\ell}(\NVE)$ & $0$ & $1$ & $1$ & $2$ & $2$ & $2$
\end{tabular}
 }
  \label{tab:ell-triangaligned}
\end{table}

\begin{table}
  \centering
  \caption{\ourfix{Sufficient $l$ for polygons with aligned edges up to 24
      edges}}
  \ourfix{ \begin{tabular}{c|cccccccccccccccccc}
$\NVE$ & $7$ & $8, 9$ & $10, 11$ & $12, 13$ & $14, 15$ & $16, 17$ & $18, 19$ & $20, 21$ & $22, 23$ & $24$
\\
$\hat{\ell}(\NVE)$ & $2$ & $3$ & $4$ & $5$ & $6$ & $7$ & $8$ & $9$ & $10$ & $11$
\\
$l$ & $1$ & $2$ & $2$ & $2$ & $3$ & $3$ & $3$ & $3$ & $4$ & $4$
\\
$\check{\ell}(\NVE)$ & $1$ & $2$ & $2$ & $2$ & $3$ & $3$ & $3$ & $3$ & $4$ & $4$
\end{tabular}
 }
  \label{tab:ell-hexagonaligned}
\end{table}

\ourfix{%
  To test numerically the coercivity of the bilinear form $\ahE{}{}$, we
  consider a set of polygons and we build for each of them the local stiffness
  matrix $A \in \mathbb{R}^{\NVE\times \NVE}$ such that
  $A_{ij} = \ahE{\varphi_i}{\varphi_j}$, involving virtual basis functions. The
  desired rank of such matrix is $\NVE-1$. In view of
  \cref{PiZeroEnhInjectivity,th:badpolyupperbound}, we define, for any
  $E\in\Mh$,
  \begin{align*}
    \hat{\ell}(\NVE) \text{ as the smallest $l$ such that }&
    2(l+1) \geq \NVE - 1 \,,
    \\
    \check{\ell}(\NVE) \text{ as the smallest $l$ such that }&
    (l+1)(l+2) \geq \NVE - 1 \,.
  \end{align*}
  Notice that \cref{PiZeroEnhInjectivity,th:badpolyupperbound} imply that the
  minimum $l$ that is sufficient to obtain local coercivity on $E$ satisfies
  $\check{\ell}(\NVE) \leq l \leq \hat{\ell}(\NVE)$, being
  $0\leq\dim\BadPoly{l}{E}\leq l(l+1)$.} \ourfix{ In the following, we compute
  numerically the minimum $l$ that induces the coercivity of the stiffness
  matrix for several sequences of polygons.  In \cref{tab:ell-regularPoly} we
  display $\hat{\ell}$, $\check{\ell}$ and the minimum $l$ required to obtain
  the desired rank computed for regular polygons of $n$ vertices having vertices
  $x_i = \begin{pmatrix} \cos\left(\frac{(i-1)\pi}{n}\right) &
    \sin\left(\frac{(i-1)\pi}{n}\right) \end{pmatrix}$, $i\in\{1,\ldots,n\}$. We
  can see that for these polygons $l=\hat{\ell}(\NVE)$. This suggests that for
  regular polygons the upper bound of \cref{th:badpolyupperbound} is
  verified. On the other hand, if we consider a sequence of non-regular convex
  polygons, the results in \cref{tab:ell-nonregularPoly} suggest that we can
  take $\dim\BadPoly{l}{E}=0$. The vertices of such polygons were generated by
  sampling random points on a circle of radius $1$ and imposing that the ratio
  of each edge and the diameter of the circle is $\geq 0.15$. A third test
  considers a sequence of polygons with aligned edges obtained starting from a
  non-equilateral triangle and then progressively splitting its edges into equal
  parts one at a time until all three edges are split into three equal parts. In
  \cref{tab:ell-triangaligned} we can see how the sufficient $l$ that guarantees
  coercivity in this case is inside the range
  $[\check{\ell}(\NVE),\hat{\ell}(\NVE)]$. A similar test is reported in
  \cref{tab:ell-hexagonaligned}, where the same procedure has been applied to a
  non-regular hexagon, thus generating a sequence of polygons up to $24$
  edges. We can see that in this case $\check{\ell}(\NVE)$ is
  sufficient. Finally, we consider a sequence of polygons that are non
  convex. To generate this sequence, we start from the quadrilateral considered
  in the second test (\cref{tab:ell-nonregularPoly}), add the edge midpoints as
  vertices and move them towards its barycenter $x_C$ with the transformation
  $S(x) = (1-\alpha)x + \alpha x_C$, thus obtaining a sequence of non-convex
  octagons. We select four polygons by choosing $\alpha
  \in\{0,0.2,0.4,0.6\}$. In all these cases, the sufficient $l$ that guarantees
  coercivity is $\check{\ell}(8) = 2$. The coordinates of all polygons
  considered in this section, except for the regular ones, are provided as
  supplementary materials to the paper.}

\subsection{Convergence tests}
\label{sec:numres:conv-test}
Let us consider problem \cref{eq:modelProblem} on the unit square with
homogeneous Dirichlet boundary conditions and the right-hand side defined such
that the exact solution is
\begin{equation*}
  U_{ex} = \sin (2\pi x) \sin (2\pi y).
\end{equation*}
In the following, we show, in log-log scale plots, the convergence curves of the
$\mathrm{L}^2$ and $\sobh{1}{}$ errors that we measure respectively as follows,
\begin{equation*}
  \mathrm{L}^2 \mbox{ error } = \sqrt{\sum_{E\in\Mh}\limits\norm[\lebl{E}]{\proj[\nabla]{1}{E}u-U_{ex}}^2} ,
\end{equation*}
\begin{equation*}
  \sobh{1}{} \mbox{ error } = \sqrt{\sum_{E\in\Mh}\limits\norm[\lebl{E}]{\nabla\proj[\nabla]{1}{E} u-\nabla U_{ex}}^2} ,
\end{equation*}
where $u$ is the discrete solution of \cref{eq:discrVarForm}.
Then, for each polygon $E\in\Mh$ we choose $l$ such that the sufficient
condition \cref{eq:suffCondWellPos} is satisfied\ourfix{, as detailed below.}

\subsubsection{Meshes}
\begin{figure}[ht]
  \centering
  \subfloat[\textit{Hexagonal}]{\label{subfig:HexagonalMesh}\includegraphics[
    scale=.6 ]{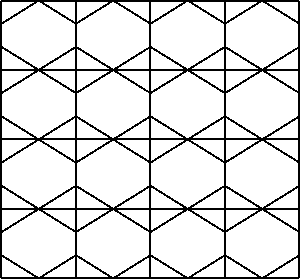} } \hspace{.05\linewidth}
  \subfloat[\textit{Octagonal}]{\label{subfig:OctagonalMesh}\includegraphics[
    scale=.6
    ]{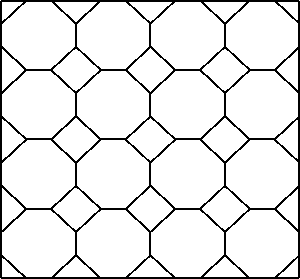}} \\
  \subfloat[\textit{Hexadecagonal}]{\label{subfig:HexadecagonalMesh}\includegraphics[
    scale=.6 ]{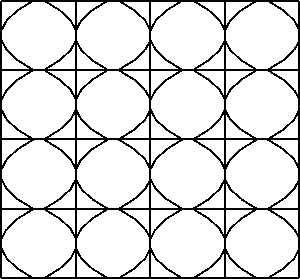}} \hspace{.05\linewidth} \subfloat[\textit{Star
    Concave}]{\label{subfig:StarConcaveMesh}\includegraphics[
    scale=.6 ]{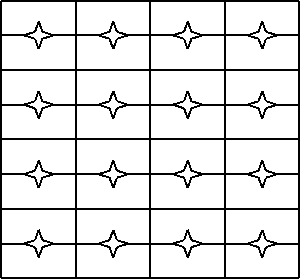}}
  \caption{Meshes}
\end{figure}

We consider four sequences of meshes for the convergence test.  The first
sequence, labeled \textit{Hexagonal}, is a tesselation made by hexagons and
triangles, as it is shown in \cref{subfig:HexagonalMesh}. \ourfix{For this mesh,
  we choose $l=0$ on triangles and $l = \hat{\ell}(6) = 2$ on hexagons.} The
second sequence, shown in \cref{subfig:OctagonalMesh} and labeled
\textit{Octagonal}, is made by octagons, squares and triangles. \ourfix{We
  choose $l=0$ on triangles, $l=\check{\ell}(4) = 1$ on squares,
  $l=\check{\ell}(8) = 2$ on octagons.} Then, the third sequence, labeled
\textit{Hexadecagonal}, is made by hexadecagons and concave pentagons, as it is
shown in \cref{subfig:HexadecagonalMesh}. \ourfix{We choose
  $l=\check{\ell}(5)=1$ on the concave pentagons and $l=\check{\ell}(16)=3$ on
  hexadecagons.} Finally, the last sequence, labeled \textit{Star Concave}, is a
non-convex tessellation made by octagons and nonagons, as it is shown in
\cref{subfig:StarConcaveMesh}. \ourfix{Here we choose $l=\hat{\ell}(8)=3$ on
  octagons and $l=\check{\ell}(9)=2$ on nonagons.} \ourfix{The choices of $l$
  were done based on a numerical evaluation of the rank of local stiffness
  matrices, as done in \cref{sec:numres:coercivity}.}
In each case we start from a mesh of $\# \Mh$ polygons then we refine it,
obtaining meshes made by $4\# \Mh$, $16\# \Mh$ and $64\# \Mh$ polygons. The
first and the third sequence start with $\# \Mh$ equal to $320$, the second and
the fourth with $\# \Mh$ equal to $164$ and $192$ respectively.

\subsubsection{Convergence results}
\begin{figure}[ht]
  \centering \subfloat[\textit{$\sobh{1}{}$
    error}]{\label{subfig:H1error}\includegraphics[
    scale=.25 ]{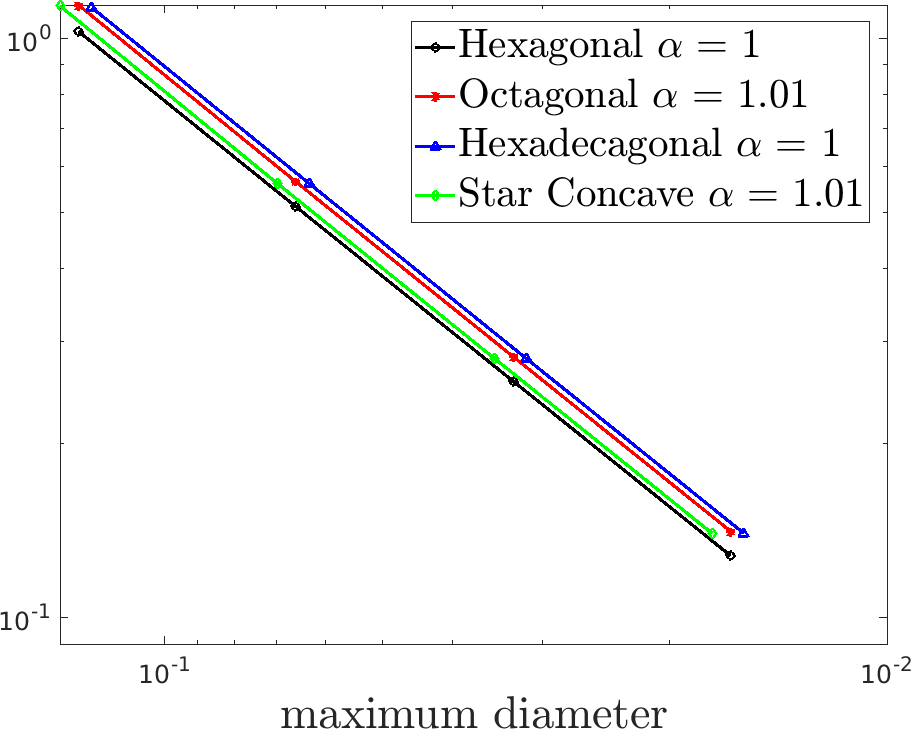} } \subfloat[\textit{$\lebl{}$
    error}]{\label{subfig:L2error}\includegraphics[
    scale=.25 ]{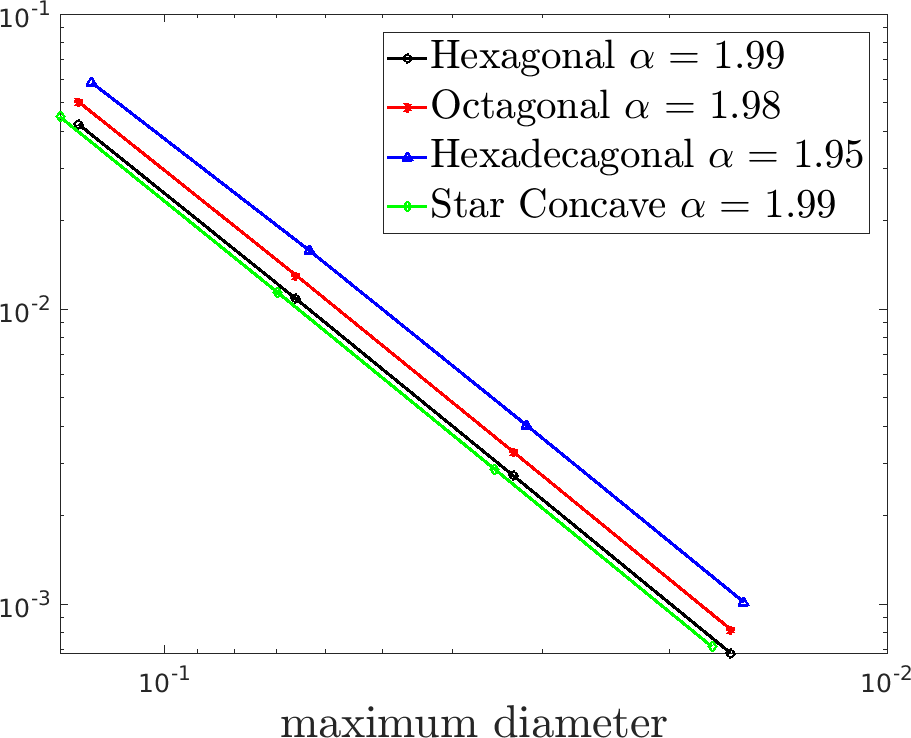}}
  \caption{Logarithmic convergence plots}
\end{figure}

For the four mesh sequences, we report the trend of the $\sobh{1}{}$ and the
$\mathrm{L}^2$ errors in \cref{subfig:H1error} and in \cref{subfig:L2error},
respectively, decreasing the maximum diameter of the polygons. In the legends,
we report the computed convergence rates with respect to $h$, denoted by
$\alpha$. We see that we get the expected values for all the meshes, as obtained
in \cref{eq:H1estimate} and \cref{eq:L2estimate}.

\fix{%
  \subsubsection{Convergence of diffusion-reaction discrete problem}
  \label{sec:numres:diffreact}
  \begin{figure}
    \centering%
    \subfloat[\textit{$\sobh{1}{}$ error}]{\includegraphics[
      scale=.25 ]{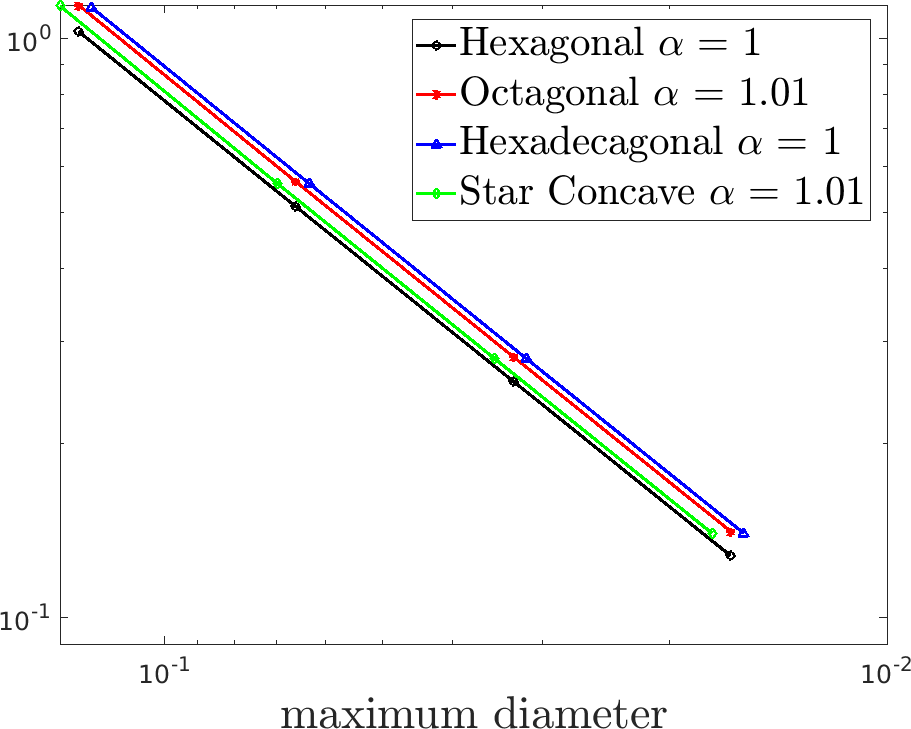} }%
    \subfloat[\textit{$\lebl{}$ error}]{\includegraphics[
      scale=.25 ]{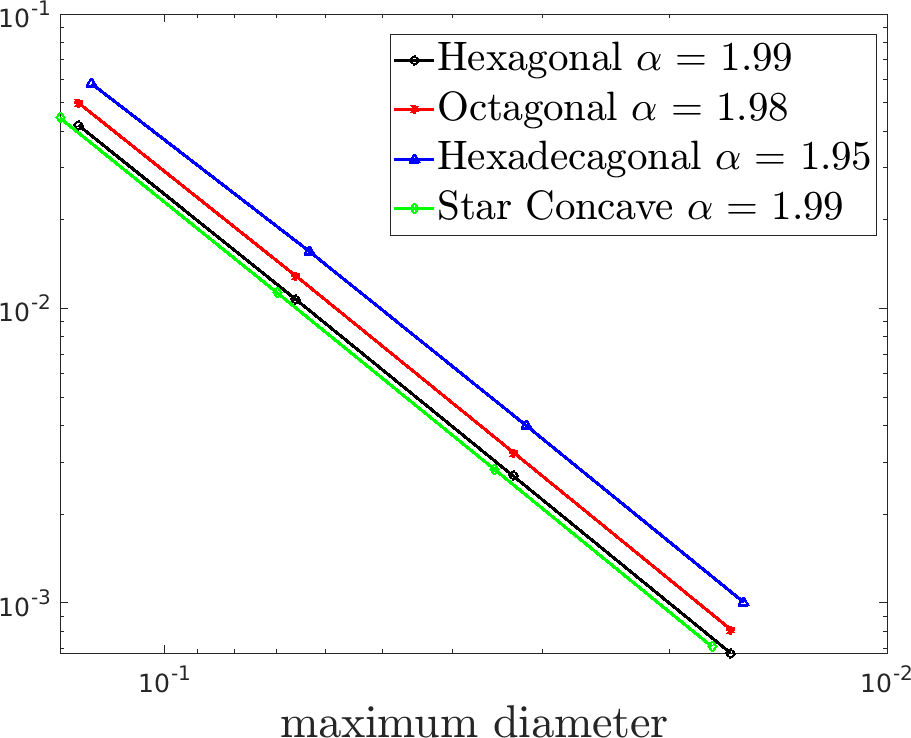}}
    \caption{\fix{Logarithmic convergence plots for diffusion-reaction model}}
    \label{fig:error-diffreact}
  \end{figure}
  We finally report, in \cref{fig:error-diffreact}, the $H^1$ and $L^2$ errors
  obtained for the four mesh sequences when solving \cref{eq:diffreact} using
  the discrete formulation \cref{eq:discr-diffreact}. We can see that the
  convergence rates $\alpha$ reported in the legends are optimal.}


\bibliographystyle{plain}
\bibliography{bibliografia}

\appendix

\section{Supplementary materials}
\subsection{Proof of \cref{lem: continuity of b}}
In order to show the proof, we have to present a preliminary result.
\begin{lemma} \label{q-DOF scaling on E}

  Let $\bar{q}\in\RQE$. Then $\exists C>0$, independent of $h_E$, such that
  \begin{equation}\label{relation of lemma q-DOF scaling on E}
    \sum_{i=1}^{\NVE} \abs{\bar{q}(x_i)} \leq C \sqrt{\sum_{\tau\in\Triang{E}}\norm[\lebl{\tau}]{\nabla \bar{q}}^2} \,.
  \end{equation}
\end{lemma}
\begin{proof}
  We notice that
  \begin{equation}
    \label{eq:q-DOF scaling on E:first-step}
    \sum_{i=1}^{\NVE} \abs{\bar{q}(x_i)} = \frac{1}{2}\sum_{\tau\in\Triang{E}}
    \left(
      \abs{\bar{q}(x_{\tau,1})} + \abs{\bar{q}(x_{\tau,2})}
    \right) \,,
  \end{equation}
  where $x_{\tau,1}$ and $x_{\tau,2}$ are the vertices of $\tau$ that are on $\partial E$. We have
  that
  \begin{equation*}
    \eval{\bar{q}}{\tau} \in \tilde{\mathbb{P}}_1(\tau) =
    \left\{p\in\Poly{1}{\tau}\colon p(x_C)=0\right\} \,,
  \end{equation*}
  and
  \begin{equation*}
    \abs{\bar{q}(x_{\tau,1})} + \abs{\bar{q}(x_{\tau,2})} =
    \norm[\ourfix{l^1}]{\dof{\tilde{\mathbb{P}}_1(\tau)}{\eval{\bar{q}}{\tau}}} \,,
  \end{equation*}
  having chosen the values at $x_{\tau,1}$ and $x_{\tau,2}$ as set of degrees of freedom on
  $\tilde{\mathbb{P}}_1(\tau)$ and denoting by $\dof{\tilde{\mathbb{P}}_1(\tau)}{\cdot}$ the
  operator returning the vector of such values. Using the mapping \cref{def:affinemapping} we get
  \begin{equation*}
    \norm[\ourfix{l^1}]{\dof{\tilde{\mathbb{P}}_1(\tau)}{\eval{\bar{q}}{\tau}}} =
    \norm[\ourfix{l^1}]{\dof{\tilde{\mathbb{P}}_1(\hat{\tau})}{\eval{\hat{\bar{q}}}{\hat{\tau}}}} \,.
  \end{equation*}
  The right-hand side of the above equation is a norm on $\tilde{\mathbb{P}}_1(\hat{\tau})$, as well
  as $\norm[\lebl{\hat{\tau}}]{\hat{\nabla}\hat{\bar{q}} }$. Then, by standard arguments about the
  equivalence of norms in finite dimensional spaces, we have
  \begin{equation*}
    \norm[\ourfix{l^1}]{\dof{\tilde{\mathbb{P}}_1(\hat{\tau})}{\eval{\hat{\bar{q}}}{\hat{\tau}}}} \leq
    \frac{\sqrt{2} \max_{i=1,2} \norm[\ourfix{l^1}]{\dof{\tilde{\mathbb{P}}_1(\hat{\tau})}{\hat{\chi}_i}}}
    {\min_{ \hat{w}\in \tilde{\mathbb{P}}_1(\hat{\tau})\colon
        \hat{w}(\hat{x}_{\hat{\tau},1})^2 + \hat{w}(\hat{x}_{\hat{\tau},2})^2=1 }
      \norm[\lebl{\hat{\tau}}]{\hat{\nabla}\hat{w}}}
    \norm[\lebl{\hat{\tau}}]{\hat{\nabla}\hat{\bar{q}} } \,,
  \end{equation*}
  where the $\hat{\chi}_i$ are Lagrangian in the degrees of freedom. Then,
  $ \norm[\ourfix{l^1}]{\dof{\tilde{\mathbb{P}}_1(\hat{\tau})}{\hat{\chi}_1}}=
  \norm[\ourfix{l^1}]{\dof{\tilde{\mathbb{P}}_1(\hat{\tau})}{\hat{\chi}_2}} = 1 $ and
  \begin{equation*}
    \norm[\ourfix{l^1}]{\dof{\tilde{\mathbb{P}}_1(\hat{\tau})}{\eval{\hat{\bar{q}}}{\hat{\tau}}}}
    \leq \frac{\sqrt{2} }{\min_{ \hat{w}\in \tilde{\mathbb{P}}_1(\hat{\tau}) \colon
        \hat{w}(\hat{x}_{\hat{\tau},1})^2 + \hat{w}(\hat{x}_{\hat{\tau},2})^2=1 }
      \norm[\lebl{\hat{\tau}}]{\hat{\nabla}\hat{w}}}
    \norm[\lebl{\hat{\tau}}]{\hat{\nabla}\hat{\bar{q}} } \,.
  \end{equation*}
  It can be proved by standard arguments that the constant in the above inequality is continuous
  with respect to $\hat{\tau}$, since it depends continuously on the deformation of the domain (see
  the proofs of \cite[Lemma 4.9]{Cangiani2016} and \cite[Lemma 4.5]{Beirao2015}). It follows by
  compactness of the set of admissible reference elements, denoted by $\Sigma$,
  (\cref{lem:compactnessSigma}) that there exists $M>0$ such that
  \begin{equation*}
    M = \max_{\hat{\tau}\in\Sigma} \frac{\sqrt{2} }
    {
      \min_{ \hat{w}\in \tilde{\mathbb{P}}_1(\hat{\tau})\colon
        \hat{w}(\hat{x}_{\hat{\tau},1})^2 + \hat{w}(\hat{x}_{\hat{\tau},2})^2=1 }
      \norm[\lebl{\hat{\tau}}]{\hat{\nabla}\hat{w}}
    } \,,
  \end{equation*}
  and thus, starting again from \cref{eq:q-DOF scaling on E:first-step} and applying the mapping
  \cref{def:affinemapping}, we get
  \begin{equation*}
    \begin{split}
      \sum_{i=1}^{\NVE} \abs{\bar{q}(x_i)} &= \frac{1}{2}\sum_{\tau\in\Triang{E}}
      \norm[\ourfix{l^1}]{\dof{\tilde{\mathbb{P}}_1(\tau)}{\eval{\bar{q}}{\tau}}} =
      \frac{1}{2}\sum_{\hat{\tau}\in\Triang{\hat{E}}}
      \norm[\ourfix{l^1}]{\dof{\tilde{\mathbb{P}}_1(\hat{\tau})}{\eval{\hat{\bar{q}}}{\hat{\tau}}}}
      \\
      &\leq \frac{M}{2}\sum_{\hat{\tau}\in\Triang{\hat{E}}}
      \norm[\lebl{\hat{\tau}}]{\hat{\nabla}\hat{\bar{q}} } = \frac{M}{2}\sum_{\tau\in\Triang{E}}
      \norm[\lebl{\tau}]{\nabla\bar{q} }
      \\
      &\leq \frac{M\sqrt{\NVE}}{2} \sqrt{\sum_{\tau\in\Triang{E}}\norm[\lebl{\tau}]{\nabla
          \bar{q}}^2} \,,
    \end{split}
  \end{equation*}
  and we obtain \cref{relation of lemma q-DOF scaling on E} since $\NVE$ is uniformly bounded by
  \cref{eq:numVerticesBounded}.
\end{proof}


Now, we can present the proof of \cref{lem: continuity of b}.
\begin{proof}
Let $\bar{q}\in\RQE$ and $\boldsymbol{v}\in\ourfixBis{\spaceV}$ be given.
Starting from \cref{eq:BAfterDivTheorem} and applying the triangular inequality, we have
\begin{equation} \label{splitting della continuita di b}
\abs{b(\bar{q},\boldsymbol v)} \leq \abs{\sum_{\tau\in\Triang{E}} \limits \int_{ \tau} \left[\nabla\bar{q}\,\boldsymbol{v} + \bar{q}\,\div \boldsymbol{v} \right] \, dx} + \abs{\sum_{i=1}^{\NVE} \limits \int_{e_i} \limits \trace{\bar{q}}{e_i}{}\jmp[e_i]{\boldsymbol{v}}\cdot \boldsymbol{n}^{e_i} ds}.
\end{equation}
Let us consider separately the two terms involved in the inequality. The first part can be analysed applying the property,
\begin{equation*}
\forall \bar{q}\in\RQE, \quad
\sum_{\tau\in\Triang{E}} \left(\norm[\lebl{\tau}]{\bar{q}} + \norm[\lebldouble{\tau}]{\nabla \bar{q}}\right)\leq\sqrt{2\NVE}\norm[\HT{E}]{\bar{q}}
\end{equation*}
and the mesh assumption \cref{eq:numVerticesBounded}, as follows
\begin{align*}
  \abs{\sum_{\tau\in\Triang{E}} \limits \int_{ \tau} \left[\nabla\bar{q}\,\boldsymbol{v} + \bar{q}\,\div \boldsymbol{v} \right] \, dx}
  &\leq \sum_{\tau\in\Triang{E}} \limits
    \left(
    \norm[\lebldouble{\tau}]{\nabla\bar{q}}\norm[\lebldouble{\tau}]{\boldsymbol{v}}+
    \norm[\lebl{\tau}]{\bar{q}}\norm[\lebl{\tau}]{\div \boldsymbol{v}}
    \right)
  \\
  &\leq C\sum_{\tau\in\Triang{E}} \limits \left(\norm[\lebldouble{\tau}]{\boldsymbol{v}}+\norm[\lebl{\tau}]{\ourfixBis{\div}\boldsymbol{v}}\right)
  \\
  &\quad \times\left(\norm[\lebldouble{\tau}]{\nabla\bar{q}}+ \norm[\lebl{\tau}]{\bar{q}}\right)
  \\
  &\leq C \norm[\HT{E}]{\bar{q}}\sum_{\tau\in\Triang{E}}\limits \left(\norm[\lebldouble{\tau}]{\boldsymbol{v}}+\norm[\lebl{\tau}]{\ourfixBis{\div}\boldsymbol{v}}\right) .
\end{align*}
Moreover, let us consider the second term of \cref{splitting della continuita di b},  computing exactly the term $\norm[\lebl{e_i}]{\trace{\bar{q}}{e_i}{}}$ and applying the properties $\forall \boldsymbol{v}\in\ourfixBis{\spaceV}$
\begin{align*}
&\sum_{i=1}^{\NVE}\limits \norm[\lebl{e_i}]{ \jmp[e_i]{\boldsymbol{v}}} \leq \sqrt{2\NVE}\sqrt{\sum_{i=1}^{\NVE}\limits \norm[\lebl{e_i}]{ \jmp[e_i]{\boldsymbol{v}}}^2}\, , \\
 &  \norm[\lebl{e_i}]{\jmp[e_i]{\boldsymbol{v}}}^2 \leq h_E \norm[\leblinfI]{\jmp[\intEdges]{\boldsymbol{v}}}^2, \; \forall e_i\in\intEdges \,,
\end{align*}
we have
\begin{align*}
\abs{\sum_{i=1}^{\NVE} \limits \int_{e_i} \limits \trace{\bar{q}}{e_i}{}\jmp[e_i]{\boldsymbol{v}}\cdot \boldsymbol{n}^{e_i} ds}
&\leq \sum_{i=1}^{\NVE}\limits \norm[\lebl{e_i}]{\trace{\bar{q}}{e_i}{}} \norm[\lebl{e_i}]{\jmp[e_i]{\boldsymbol{v}}\cdot \boldsymbol{n}^{e_i}}  \\
&\leq \sum_{i=1}^{\NVE}\limits \frac{\sqrt{h_{e_i}}}{\sqrt{3}}\abs{\bar{q}(x_i)}\norm[\lebldouble{e_i}]{\jmp[e_i]{\boldsymbol{v}}} \\
&\leq \frac{h_E}{\sqrt{3}} \norm[\leblinfI]{\jmp[\intEdges]{\boldsymbol{v}}} \sum_{i=1}^{\NVE}\limits \abs{\bar{q}(x_i)} \\
& \leq C h_E \norm[\leblinfI]{\jmp[\intEdges]{\boldsymbol{v}}} \norm[\HT{E}]{\bar{q}}\, ,
\end{align*}
where we apply \cref{q-DOF scaling on E} in the last step.
Finally, substituting
into \cref{splitting della continuita di b}, we obtain
\begin{align*}
  \abs{b(\bar{q},\boldsymbol v)}
  &\leq C \norm[\HT{E}]{\bar{q}} \left(
    \sum_{\tau\in\Triang{E}}\limits \left(\norm[\lebldouble{\tau}]{\boldsymbol{v}} +
    \norm[\lebl{\tau}]{\ourfixBis{\div}\boldsymbol{v}}\right) +
    h_E  \norm[\leblinfI]{\jmp[\intEdges]{\boldsymbol{v}}}\right)
  \\
  &\leq C\norm[\HT{E}]{\bar{q}} \norm[\ourfixBis{\spaceV}]{\boldsymbol{v}}\,.
\end{align*}
\end{proof}


\subsection{Proof of \cref{thm:L2estimate}}
\begin{proof}
Let us define the auxiliary problem: let $\Psi\in\sobh{2}{\Omega}\cap\sobh[0]{1}{\Omega}$ the solution of $\a{V}{\Psi}=\scal[\Omega]{U-u}{V}$ $\forall V\in\sobh[0]{1}{\Omega}$. From the definition of $\Psi$, we get: 
\begin{align}\label{eq:StimaPsi}
\exists C>0 : \quad \seminorm[2]{\Psi}&\leq C\norm[\lebl{\Omega}]{U-u} ,
\\
\label{eq:StimaPsi-1}
\exists C>0 : \quad \norm[{\sobh[0]{1}{\Omega}}]{\Psi}&\leq C\norm[\lebl{\Omega}]{U-u}.
\end{align}
Let us denote by $\Psi_I$ the interpolant of $\Psi$ according to \cref{lem:interpolationError}. 
Applying the auxiliary problem, the discrete problem \cref{eq:discrVarForm} and the definition of the bilinear form $a$ \cref{eq:defBilFormA}, we have
\begin{align}
  \notag
  \norm[\lebl{\Omega}]{U-u}^2
  &= \scal[\Omega]{U-u}{U-u} = \a{U-u}{\Psi}
  \\ \notag
  &= \a{U}{\Psi-\Psi_I}+\a{U}{\Psi_I}-\a{u}{\Psi}
  \\ \notag
  &=\a{U}{\Psi-\Psi_I}+\scal[\Omega]{f}{\Psi_I}-\a{u}{\Psi}
  \\ \notag
  &\begin{gathered}
    =\a{U}{\Psi-\Psi_I}+\scal[\Omega]{f}{\Psi_I}-\left(\sum_{E\in\Mh}\scal[E]{f}{\proj{0}{E}\Psi_I}\right)
    \\
    + \ah{u}{\Psi_I}-\a{u}{\Psi} + \a{u}{\Psi_I} - \a{u}{\Psi_I}
  \end{gathered}
  \\
  &\begin{gathered}
    = \a{U-u}{\Psi-\Psi_I}+\left(\sum_{E\in\Mh}\scal[E]{f}{\Psi_I-\proj{0}{E}\Psi_I}\right)
    \\
    +\ah{u}{\Psi_I}-\a{u}{\Psi_I} .
  \end{gathered} \label{eq:Step1Proof}
\end{align}
Let us consider the terms of the previous relation separately.
First, applying the Cauchy-Schwarz inequality, \cref{eq:interpolationError}, \cref{eq:polynomialErrorPiZero} and \cref{eq:StimaPsi}, we have, for the first term,
\begin{align} \notag
\a{U-u}{\Psi-\Psi_I}&\leq \norm[{\sobh[0]{1}{\Omega}}]{U-u}\norm[{\sobh[0]{1}{\Omega}}]{\Psi-\Psi_I} \\ &\leq Ch\norm[{\sobh[0]{1}{\Omega}}]{U-u}\seminorm[2]{\Psi}
\leq Ch\norm[{\sobh[0]{1}{\Omega}}]{U-u}\norm[\lebl{\Omega}]{U-u} ,
\label{eq:firstTerm}
\end{align}
and, for the second one,
\begin{align}\notag
\sum_{E\in\Mh}\scal[E]{f}{\Psi_I-\proj{0}{E}\Psi_I}&=\sum_{E\in\Mh}\scal[E]{f-\proj{0}{E} f}{\Psi_I-\proj{0}{E}\Psi_I}\\ \notag
&\leq\sum_{E\in\Mh}\norm[\lebl{E}]{f-\proj{0}{E} f}\norm[\lebl{E}]{\Psi_I-\proj{0}{E}\Psi_I}\\
&\leq Ch\seminorm[\sobh{1}{\Omega}]{f}\sum_{E\in\Mh}\norm[\lebl{E}]{\Psi_I-\proj{0}{E}\Psi_I}. \label{eq:Step1-forceTerm}
\end{align}
Applying the property 
\begin{equation*}
\forall E\in\Mh, \;\;
\norm[\lebl{E}]{\Psi_I-\proj{0}{E}\Psi_I}\leq \norm[\lebl{E}]{\Psi_I-\proj{0}{E}\Psi},
\end{equation*}
 \cref{eq:interpolationError} and \cref{eq:polynomialErrorPiZero} to \cref{eq:Step1-forceTerm}, we obtain
\begin{align}\notag
\sum_{E\in\Mh}\scal[E]{f}{\Psi_I-\proj{0}{E}\Psi_I}&\leq Ch\seminorm[\sobh{1}{\Omega}]{f}\sum_{E\in\Mh}\norm[\lebl{E}]{\Psi_I-\proj{0}{E}\Psi}\\ \notag
&\leq Ch\seminorm[\sobh{1}{\Omega}]{f}\sum_{E\in\Mh}\left(\norm[\lebl{E}]{\Psi_I-\Psi}+\norm[\lebl{E}]{\Psi-\proj{0}{E}\Psi}\right) \\
&\leq Ch\seminorm[\sobh{1}{\Omega}]{f}\left(h^2\seminorm[2]{\Psi} + h\norm[{\sobh[0]{1}{\Omega}}]{\Psi}\right).\label{eq:Step2-forceTerm}
\end{align}
We can omit higher order terms and apply \cref{eq:StimaPsi-1}, obtaining
\begin{equation}\label{eq:SecondTerm}
\sum_{E\in\Mh}\scal[E]{f}{\Psi_I-\proj{0}{E}\Psi_I}\leq Ch^2\seminorm[\sobh{1}{\Omega}]{f}\norm[\lebl{\Omega}]{U-u}.
\end{equation}
Finally, we have to bound $\ah{u}{\Psi_I}-\a{u}{\Psi_I}$. Then, applying the orthogonality property of $\proj{l}{E}$, adding and subtracting terms,
we have
\begin{align}
  \notag \ah{u}{\Psi_I}-\a{u}{\Psi_I}
  &= \sum_{E\in\Mh}\scal[E]{\proj{l}{E}\nabla
    u}{\nabla \Psi_I}-\scal[E]{\nabla u}{\nabla\Psi_I}
  \\
  \notag
  &=\sum_{E\in\Mh}\scal[E]{\proj{l}{E}\nabla u-\nabla u}{\nabla
    \Psi_I-\proj{\ourfix{0}}{E}\nabla\Psi_I}
  \\
  &\begin{gathered} =
    \sum_{E\in\Mh}\scal[E]{\proj{l}{E}\nabla
      u-\proj{l}{E}\nabla U}{\nabla
      \Psi_I-\proj{\ourfix{0}}{E}\nabla\Psi_I} \\ +
    \scal[E]{\proj{l}{E}\nabla U-\nabla U}{\nabla
      \Psi_I-\proj{\ourfix{0}}{E}\nabla\Psi_I}
    \\
    + \scal[E]{\nabla U-\nabla u}{\nabla
      \Psi_I-\proj{\ourfix{0}}{E}\nabla\Psi_I}. \label{eq:Step1:lastTerm}
  \end{gathered}
\end{align}
Notice that, applying \cref{eq:interpolationError} and \cref{eq:polynomialErrorGrad}, we have the property $\forall E\in\Mh$ :
\begin{equation*}
\norm[\lebl{E}]{\nabla \Psi_I-\proj{\ourfix{0}}{E}\nabla\Psi_I} \leq \norm[\lebl{E}]{\nabla \Psi_I-\proj{\ourfix{0}}{E}\nabla\Psi} 
\leq Ch\seminorm[2,E]{\Psi} .
\end{equation*}
Therefore, applying the continuity of the projection operator and \cref{eq:StimaPsi},  the first and the last term of \cref{eq:Step1:lastTerm} can be bounded as
\begin{equation} \label{eq:LastTerm-1}
\begin{gathered} 
\sum_{E\in\Mh}
\scal[E]{\proj{l}{E}\nabla u-\proj{l}{E}\nabla U}{\nabla \Psi_I-\proj{\ourfix{0}}{E}\nabla\Psi_I} +
\scal[E]{\nabla U-\nabla u}{\nabla \Psi_I-\proj{\ourfix{0}}{E}\nabla\Psi_I}
\\ \leq Ch\norm[{\sobh[0]{1}{\Omega}}]{U-u}\norm[\lebl{\Omega}]{U-u}.
\end{gathered}
\end{equation}
Similarly, the second term is bounded as
\begin{equation} \label{eq:LastTerm-2}
\sum_{E\in\Mh}\scal[E]{\proj{l}{E}\nabla U-\nabla U}{\nabla \Psi_I-\proj{\ourfix{0}}{E}\nabla\Psi_I} \leq C h^2 \seminorm[2]{U}\norm[\lebl{\Omega}]{U-u} .
\end{equation}
Finally, applying \cref{eq:firstTerm},\cref{eq:SecondTerm},\cref{eq:LastTerm-1} and \cref{eq:LastTerm-2} to \cref{eq:Step1Proof} and simplifying, we obtain
\begin{equation*}
\norm[\lebl{\Omega}]{U-u} \leq C \left( h\norm[{\sobh[0]{1}{\Omega}}]{U-u} +  h^2\seminorm[\sobh{1}{\Omega}]{f} +  h^2 \seminorm[2]{U}\right) .
\end{equation*}
Applying the $\sobh{1}{}$-estimate (\cref{thm:H1estimate}) we obtain the relation \cref{eq:L2estimate}.
\end{proof}

\end{document}